\documentclass[journal,twoside,web]{IEEEtran}

    \usepackage[pdftex]{graphicx}
    \DeclareGraphicsExtensions{.pdf,.jpeg,.png}

\usepackage{generic}
\usepackage{epstopdf}
\usepackage{graphicx}
\usepackage{amssymb}
\usepackage{amsmath}
\usepackage{float}
\usepackage{amsfonts, bm}
\usepackage{color}
\usepackage[noadjust]{cite}

\newcommand{\Lp}{\mathbb{L}_p}
\newcommand{\Leb}{\mathbb{L}}

\newcommand{\Hp}{\mathbb{H}_p}
\newcommand{\Htwo}{\mathbb{H}_2}
\newcommand{\Hinf}{\mathbb{H}_\infty}
\newcommand{\Reals}{\mathbb{R}}
\newcommand{\Complexes}{\mathbb{C}}

\newcommand{\map}[1]{\bm{#1}}

\renewcommand{\i}{u}
\renewcommand{\v}{v}
\newcommand{\R}{R}

\newcommand{\tauS}{\tau_s}
\newcommand{\tauR}{\tau_r}

\newcommand{\Y}{Y}
\newcommand{\eye}{I}
\newcommand{\Z}{Z}
\newcommand{\q}{q}
\newcommand{\p}{r}
\newcommand{\G}{G}

\newcommand{\x}{x}
\newcommand{\A}{A}
\newcommand{\B}{B}
\newcommand{\C}{C}
\newcommand{\D}{D}
\renewcommand{\P}{P}
\newcommand{\X}{X}

\DeclareMathOperator{\spec}{spec}
\DeclareMathOperator{\rank}{rank}
\DeclareMathOperator{\blkdiag}{block diag}
\DeclareMathOperator{\Ex}{\mathcal{E}}

\renewcommand{\theenumi}{\alph{enumi}}

\newtheorem{assumption}{Assumption}
\newtheorem{th1}{Theorem}
\newtheorem{lm1}{Lemma}
\newtheorem{cor1}{Corollary}
\newtheorem{def1}{Definition}

\title{\LARGE\textbf{
On the Feasibility of Self-Powered Linear Feedback Control}
}

\author{Connor H. Ligeikis and Jeffrey T. Scruggs
\thanks{The first author was supported by an NSF Graduate Research Fellowship. This funding is gratefully acknowledged. Views expressed in this paper are those of the authors and do not necessarily reflect those of the National Science Foundation.}
\thanks{C. Ligeikis and J. Scruggs are with the Department of Civil \& Environmental Engineering,
      	University of Michigan, Ann Arbor, MI, 48109.
      	Phone: \mbox{734-764-1812}, email: \mbox{jscruggs@umich.edu}, \mbox{ligeikis@umich.edu}}
}

\begin{document}




\maketitle

\begin{abstract}

A control system is called \emph{self-powered} if the only energy it requires for operation is that which it absorbs from the plant. 
For a linear feedback law to be feasible for a self-powered control system, its feedback signal must be colocated with the control inputs, and its input-output mapping must satisfy an associated passivity constraint.
The imposition of such a feedback law can be viewed equivalently as the imposition of a linear passive shunt admittance at the actuation ports of the plant.
In this paper we consider the use of actively-controlled electronics to impose a self-powered linear feedback law.
To be feasible, it is insufficient that the imposed admittance be passive, because parasitic losses must additionally be overcome.
We derive sufficient feasibility conditions which explicitly account for these losses. 
In the finite-dimensional, time-invariant case, the feasibility condition distills to a more conservative version of the Positive Real Lemma, which is parametrized by various loss parameters. 
Three examples are given, in which this condition is used to determine the least-efficient loss parameters necessary to realize a desired feedback law.
\end{abstract}
\begin{IEEEkeywords}
Passivity, Energy Systems, Constrained Control, Linear Feedback
\end{IEEEkeywords}

\section{Introduction}

\begin{figure}
    \centering
    \includegraphics[scale=.9]{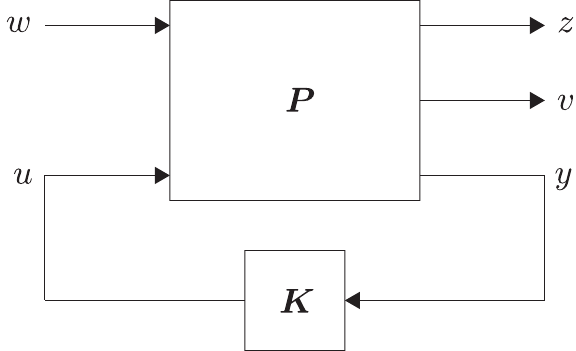}
    \caption{Generic block diagram}
    \label{block}
\end{figure}

In many applications of feedback control to physical systems, closed-loop performance objectives can be achieved merely by extracting, storing, and re-injecting the energy injected into a plant by exogenous disturbances. 
To make this statement more precise, consider Fig.~\ref{block} which shows a plant $\map{P} : \{w,u\} \mapsto \{v,y,z\}$ controlled by a feedback law $\map{K} : y \mapsto u$. 
Vector $u(t) \in \mathbb{R}^{n_p}$ is the control input, and the associated output vector $v(t) \in \mathbb{R}^{n_p}$ is such that the inner product $u^T(t) v(t)$ is the power injected into $\map{P}$ by the control inputs at time $t$.
Physically, each input/output pair $\{u_i,v_i\}$ constitutes the ``flow'' and ``effort'' variables associated with port $i \in \{1,...,n_p\}$, embedded within $\map{P}$.
For example, if $u_i(t)$ is an electric current into port $i$, then $v_i(t)$ is the corresponding voltage across the port terminals. 
Likewise, if $u_i(t)$ is the mechanical force between two degrees of freedom in a plant, then $v_i(t)$ is the relative velocity between these degrees of freedom.
We assume that the instantaneous power that must be expended to implement a control input vector $\i(t)$ is equal to $\i^T(t)R\i(t)$, for some matrix $R=R^T>0$.
Vector $w(t) \in \mathbb{R}^{n_w}$ is the exogenous disturbance to the closed-loop system, $y(t)\in\mathbb{R}^{n_y}$ is the set of measurements available for feedback, and $z(t) \in \mathbb{R}^{n_z}$ is the set of closed-loop response quantities by which performance is judged. 
We presume the feedback control system is operational for all $t \geqslant 0$, and that $\i(t) = 0$ for all $t < 0$. 
In this context, this paper is motivated by problems in which the objective is to design $\map{K}$ such that the closed-loop mapping $w \mapsto z$ is favorable in some sense, subject to the constraint that the running total energy expended by the controller must be negative for all $t > 0$, i.e., 
\begin{equation}\label{intro_1}
    \int_{0}^t \left[ u^T(\tau) v(\tau) + u^T(\tau)Ru(\tau) \right] d\tau \leqslant 0 , \ \ \forall t \in \mathbb{R}_{> 0}.
\end{equation}

The advantages of imposing constraint \eqref{intro_1} on feedback law $\map{K}$ are twofold.
The first advantage is that in the ideal case, it implies that $\map{K}$ can be realized without the need for energy to be supplied.
(By ``ideal'' we mean that the control hardware is able to store and reuse the energy extracted from the plant with perfect efficiency.)
This allows such feedback laws to operate in complete energy-autonomy, a significant design advantage in applications for which energy-efficiency is important, or for which power delivery is costly or unreliable. 

The second advantage arises in applications for which $\map{P}$ is such that $z \in \mathbb{L}_2$ for all $w,u \in \mathbb{L}_2$, and where the open-loop mapping $u\mapsto v$ is passive, i.e., where there exists $\beta : \mathbb{L}_2 \mapsto \mathbb{R}_{\geqslant 0}$ such that
\begin{equation}\label{intro_2}
    \int_{-\infty}^\infty u^T(\tau) v(\tau) d\tau \geqslant -\beta(w), \ \ \forall w , u \in \mathbb{L}_2.
\end{equation}
For many physical applications, these properties can be guaranteed to hold even if significant model uncertainty exists, because their violation would be inconsistent with the laws of thermodynamics for the plant.
In such circumstances, these properties, together with \eqref{intro_1}, imply the robust input-output stability of the closed-loop system (i.e., $z \in \mathbb{L}_2$ for all $w \in \mathbb{L}_2$), as a result of the Passivity Theorem \cite{desoer2009feedback}.

Due to these advantages, the above constrained control problem has a long history in disturbance rejection problems, such as in mechanical vibration suppression \cite{
miller1990optimal,
elliott1991power,
macmartin1991control,
jolly1997assessing,
sharp2002theoretical,
kelkar2004control,
gosavi2004modelling,
zilletti2014optimisation
}.
It has also received considerable attention in the area of robotics \cite{
ortega2002stabilization,
albu2007unified,
secchi2007control,
hatanaka2015passivity,
ott2008passivity,
chevva2015active
}.
The problem also naturally arises in other situations where both $\map{P}$ and $\map{K}$ are modeled via realizations that are port-Hamiltonian with dissipation \cite{
ortega2001putting,
van2014port,
ott2008passivity,
beattie2019robust,
rashad2020twenty,
borja2020new
}. 

\subsection{Linear, Colocated, SPR Feedback Control}

In the case where $y = v$, $\map{K}$ is called \emph{colocated}. 
For the physical case in which $u(t)$ is a vector of electrical port currents and $v(t)$ is a corresponding vector of port voltages, the control system can then be viewed as imposing an effective shunt admittance $\map{Y} = -\map{K}$ at the port terminals, i.e., 
\begin{equation} \label{i=-Yv}
u(t) = -(\map{Y}v)(t)
\end{equation}
where, in order to satisfy \eqref{intro_1} for all $v \in \Leb_{2e}^+$, it is necessary that $\map{Y}$ be an output-strictly-passive mapping \cite{brogliato2019dissipative}.
Colocated control is advantageous because in many physical systems, $v(t)$ is readily available for feedback, often without the need of additional sensing hardware beyond that which is necessary for actuation of inputs $u(t)$.
Additionally, it has the convenience that the actuators and sensors are at the same physical location. 

Another appeal of colocated control is that, in the case that $\map{Y}$ is restricted to be linear and time-invariant (LTI) and finite-dimensional, it can be physically realized via shunt networks of passive electrical or mechanical components, connected to the ports.
In the case of mechanical shunts, the network is comprised of springs, inertia, dampers, and levers \cite{smith2002synthesis}.
In the case of electrical shunts, the network is comprised of inductors, capacitors, resistors, transformers (and, possibly, gyrators) \cite{newcomb1966linear}.
It is a classical result \cite{brune1931synthesis} that an LTI, finite-dimensional admittance $\map{Y}$ can be so realized if and only if it is positive-real (PR), i.e., if its transfer function $\hat{Y}(s)$ is analytic and satisfies $\hat{Y}(s)+\hat{Y}^H(s) > 0$ for all $s \in \mathbb{C}_{>0}$.
It is also a classical result that $\map{Y}$ is PR if and only if it is passive, i.e., if \eqref{intro_1} holds with $R = 0$ for all $v \in \Leb_{2e}^+$.
Synthesis of passive electrical networks to produce a desired PR $\map{Y}$ is a central topic of classical network theory, with the most prevailing techniques being the Brune \cite{brune1931synthesis} and Darlington \cite{darlington1939synthesis} synthesis procedures and the various offshoots of these \cite{newcomb1966linear,hughes2018electrical}.
Analogous network synthesis techniques exist for mechanical admittances \cite{smith2002synthesis}.
The domain of admittances $\map{Y}$ realizable with imperfect (i.e., lossy) components is the domain of Strictly Positive Real (SPR) transfer functions \cite{taylor1974strictly}.
(A transfer function $\hat{Y}(s)$ is SPR if there exists $\epsilon>0$ such that $\hat{Y}(s-\epsilon)$ is PR.)
For LTI finite-dimensional $\map{Y}$, there exists an $R =R^T > 0$ such that \eqref{intro_1} holds for all $v \in \Leb_{2e}^+$, if $\map{Y}$ is SPR and if its transfer function at infinite frequency (denoted $D$) satisfies $D + D^T \geqslant D^TRD$. 
Consequently, when the two conditions discussed previously hold for $\map{P}$, the restriction of $\map{Y}$ to LTI, finite-dimensional, SPR mappings assures robust input-output stability in the $\mathbb{L}_2$ sense.

Many studies have proposed a fusion of optimal control techniques with network synthesis techniques, resulting in various SPR-constrained colocated feedback optimization methodologies.
Some of these methods accomplish network optimization directly; i.e., they start from a prescribed passive network structure parametrized by values of its hardware components, and then optimize those parameters under some closed-loop performance measure, such as an $\Htwo$ or $\Hinf$ objective \cite{zuo2003structured,chen2009restricted}.
Other approaches are indirect, and take advantage of network synthesis techniques.
In these techniques, the optimal SPR-constrained $\map{Y}$ is solved first, in the absence of an assumed network structure \cite{lozano1988design,
papageorgiou2006positive,
geromel1997synthesis,
damaren2006optimal,
bridgeman2014conic,
warner2015control,
forbes2019synthesis}. 
Then, the associated physical network is synthesized from the optimal admittance via a standard network synthesis techniques (i.e., Brune, Darlington, etc.).
Irrespective of whether direct  or indirect techniques are used, the resultant optimization of $\map{Y}$ is in general nonconvex.
As yet, there is no known technique for converting either approach into a convex problem, without the introduction of artificial conservatism.

\subsection{Linear, Colocated, Self-Powered Feedback Control}

It is straightforward to show that if $\map{K}$ is linear, then it must be colocated in order to satisfy \eqref{intro_1} for all $v \in \Leb_{2e}^+$.
As such, all linear controllers adhering to this constraint may be characterized via an output-strictly passive admittance $\map{Y}$.
However, there are some disadvantages to the realization of a desired $\map{Y}$ via passive shunt networks, as discussed above.
Firstly, acceptable performance may necessitate impractical network realizations. 
Additionally, passive shunt networks have the shortcoming of being static, i.e., they are designed once, and are then unable to adapt to changes.  

For these reasons, there are advantages to the realization of $\map{Y}$ synthetically (i.e., with actively-controlled electronics) even if it could theoretically be realized via passive shunts \cite{jolly1997assessing, fleming2000synthetic, moheimani2003survey, sugino2018design}.
Such an approach still exhibits stability-robustness, while also affording the advantages of time-variability and comparatively-compact hardware realizations.
However, it does not guarantee energy-autonomy, because the energy storage subsystem and the electronics exhibit parasitic losses, which must be overcome in order to realize a desired $\map{Y}$.
As such, even though $\map{Y}$ may satisfy constraint \eqref{intro_1}, its implementation may still require external power.

This issue is the motivation for this paper.
Conceptually, we are interested in the feasibility of imposing a synthetic admittance $\map{Y}$, using only the energy absorbed from the plant.
We refer to such a control technology as \emph{self-powered}, and to a feasible $\map{Y}$ as a Self-Powered Synthetic Admittance (SPSA).
Self-powered control has been investigated by numerous researchers in the area of vibration suppression \cite{%
nakano2003self,%
nakano2004combined,%
khoshnoud2015energy,%
choi2009self,%
tang2011self,%
asai2016nonlinear}.
Other terms in the vibration suppression literature for these technologies include ``regenerative'' 
\cite{%
jolly1997regenerative,%
jolly1997assessing,%
margolis2005energy,
anubi2015energy,%
clemen2016regenerative,%
liu2013regenerative,%
shen2018energy}
and ``energy recycling'' 
\cite{%
onoda2003energy,%
onoda2008performance}
control systems.
In the area of robotics, self-powered control has received significant recent attention, for use in active prosthetics and autonomous systems
\cite{laschowski2019lower,
goldfarb2003design,
richter2015framework,
carabin2017review,
khalaf2018global,
richter2014semiactive,
khalaf2019trajectory
}.
The above studies vary significantly in their modeling assumptions.
(For example, in the studies on regenerative control, constraint \eqref{intro_1} is imposed only as $t\rightarrow\infty$.)
Moreover, none develop an explicit feasibility condition that accounts for parasitic dissipation in the energy storage subsystem.
Furthermore, they all presume the controller to be LTI and finite-dimensional. 
Against this backdrop, the contribution of this paper is a feasibility criterion for SPSAs which ({\itshape i}) accounts for parasitic dissipation, and ({\itshape ii}) accommodates time-varying and infinite-dimensional controllers. 

We note that a preliminary version of the theory developed in this paper was recently published in a conference paper by the authors \cite{ligeikis2021feasibility}. 
However, in \cite{ligeikis2021feasibility} the scope was restricted to finite-dimensional, time-invariant SPSAs, and the primary contribution was related to the synthesis of a feasible $\map{Y}$ under an $\mathbb{H}_2$ performance measure. 
By contrast, the present paper extends this theory in a much broader class of linear SPSAs.
Also, the present paper demonstrates the use of this theory to characterize Pareto-optimal parasitic loss parameters to realize a given passive feedback law as an SPSA. 
Beyond these distinctions, the present paper also contains complete proofs, rather than the incomplete sketches provided in \cite{ligeikis2021feasibility}.

The paper is organized as follows.  
Section \ref{sec:sec1} presents the modeling assumptions made for the physical realization of an SPSA, and derives implicit feasibility conditions.
Section \ref{sec:sec2} then derives an explicit sufficiency condition for SPSA feasibility, first for a broad class of infinite-dimensional, time-varying linear SPSAs, and then for the special case in which the SPSA is time-invariant.
Section \ref{sec:sec3} examines this sufficient feasibility condition for the case in which an SPSA can be realized with a finite-dimensional state space, and re-casts the feasibility criterion in terms of matrix inequalities which are reminiscent of (but distinct from) the one which arises in the Positive Real Lemma. 
Also in this section, it is shown that the derived feasibility condition is necessary as well as sufficient, in two important special cases. 
Section \ref{sec:sec4} presents three examples in which the SPSA feasibility condition is used to identify a Pareto front of the least-efficient parasitic loss parameters necessary to realize a given passive admittance as an SPSA. 
Section \ref{sec:sec5} makes some brief conclusions.

\subsection{Notation}
Bold italics is used for generic mappings, while (possibly time-valued) matrices in $\Complexes^{n\times m}$ or $\Reals^{n\times m}$ are italicized.
Sets $\Reals_{>0}$ and $\Reals_{\geqslant 0}$ refer to segments $(0,\infty)$ and $[0,\infty)$, respectively.
The statements $X>0$ and $X\geqslant 0$ for a square matrix $X$ both imply that $X$ is Hermitian, and that its eigenvalues are in $\Reals_{>0}$ and $\Reals_{\geqslant 0}$ respectively.
For a matrix $X\in\Complexes^{m\times n}$, $X^T$ and $X^H$ are the transpose and complex conjugate transpose, respectively.
Norms $\| x \|_2$ and $\| x \|_\infty$ refer to the Euclidean and infinity norms of vector $x\in\Complexes^{n}$. 
The notation $x^Hx$ and $x^Tx$ are used interchangeably with $\|x\|_2^2$ for $x\in\Complexes^n$ and $\Reals^n$ respectively. 
Matrix norm $\| X \|_2$ refers to the maximum singular value of $X\in\Complexes^{n\times m}$. 
For a matrix function $X : \Reals \rightarrow \Complexes^{n\times m}$, the notation $\hat{X}(s)$ refers to its Laplace transform, assuming it exists.
The notation $\Lp(\mathbb{D})$ where $\mathbb{D}\subseteq \Reals^m$ and $p\in[1,\infty]$, refers to the Lebesgue space of functions $U:\mathbb{D}\rightarrow\mathbb{C}^n$ for which $\| U \|_{\Lp(\mathbb{D})} \triangleq \left\{ \int_{x\in\mathbb{D}} \|U(x)\|^p d\mathbb{D} \right\}^{1/p} < \infty$. 
The notation $\mathbb{L}_{pe}(\mathbb{D})$ refers to the extended $\Lp(\mathbb{D})$ space, for which $\left\{ \int_{x\in\mathbb{D}'} \|U(x)\|^p d\mathbb{D}' \right\}^{1/p} < \infty$ for any compact subdomain $\mathbb{D}'\subset\mathbb{D}$. 
Because of their frequent usage, we use the short-hand $\Leb_p^+$ and $\Leb_{pe}^+$ to refer to $\Leb_{p}(\Reals_{\geqslant 0})$ and $\Leb_{pe}(\Reals_{\geqslant 0})$, respectively.
For a linear, time-invariant system with transfer function matrix $\hat{X}(s)$, norm $\| \hat{X} \|_{\Hp}$ refers to its $\Hp$ norm.

\section{Modeling and Assumptions}\label{sec:sec1}

\subsection{Physical Modeling}

\begin{figure}
	\centering
	\includegraphics[scale=.9]{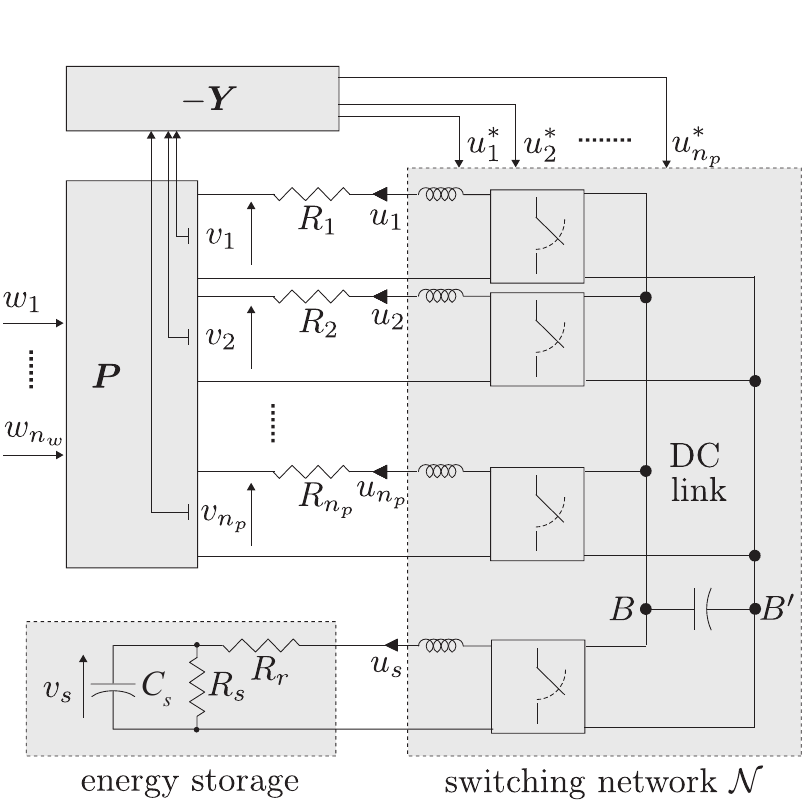}
	\caption{Embedment of an SPSA within a passive plant}
	\label{fig2}
\end{figure}

Fig.~\ref{fig2} shows a schematic for one possible implementation of a multiport SPSA, embedded within a passive plant subjected to exogenous disturbance $w(t) \in \Reals^{n_w}$.
The plant has $n_p$ electric ports, with the dynamics of each port being characterized by a voltage $v_k$, and a colocated current $u_k$, for $k \in \{1,...,n_p\}$. 
The hardware that realizes the SPSA is comprised of passive components (i.e., inductors, resistors, and capacitors) and transistors.
The transistors are operated as switches that open and close different circuit paths at high-frequency, in a pulse-width modulation (PWM) paradigm, in such a way as to control the currents $u_k$, $k \in \{1,...,n_p\}$. 
This switching network $\mathcal{N}$ interfaces each of the $n_p$ ports with a common voltage bus $B-B'$, called a \emph{DC Link}. 
The DC link allows for energy removed from the plant at one port to be re-injected into the plant at another port.
Also interfaced to the DC link is an energy-storage subsystem, depicted as capacitor $C_s$. 
(Supercapacitors, flywheels, or batteries could also be used.) 
This energy-storage subsystem is non-ideal, in the sense that it dissipates energy whenever it delivers or accepts power, and because it exhibits leakage. 
These effects are represented by resistances $R_r$ and $R_s$, respectively. 

The switching network $\mathcal{N}$ controls the current vector $\i(t)$ so as to track a desired current command vector, $\i^*(t) = -(\map{Y}v)(t)$. 
As such, by facilitating current tracking at high-bandwidth, $\mathcal{N}$ approximately imposes synthetic admittance relationship \eqref{i=-Yv}.

Three assumptions are made, regarding $\mathcal{N}$.  
These assumptions constitute a simplified view of a physical system which in reality exhibits much more nuance. 

\begin{assumption}\label{lossless_assumption}
$\mathcal{N}$ is lossless.  
\end{assumption}

As such, the dissipative electrical components in the SPSA hardware consist of resistances $R_1, R_2, ...R_{n_p}, R_r,$ and $R_s$.  
The total power dissipation, denoted $P_d(t)$, is thus
\begin{equation} \label{eq0_1}
	P_d(t) = \i^T(t)\R\i(t) + R_r u_s^2(t) + \tfrac{1}{R_s}v_s^2(t)
\end{equation}
where $\R \triangleq \text{diag}\left\{R_1...R_{n_p}\right\}$.
We note that Assumption \ref{lossless_assumption} does not require that the physical realization of the power-electronics be lossless. 
This is because for modeling purposes, conductive losses in the switching blocks may be subsumed into resistances $R_1,...,R_{n_p}$, as well as $R_r$.
As such, in reality these resistances account for the dissipation of the power-electronics, as well as the actuators and energy storage subsystem.
Nonetheless, it is important to note that Assumption \ref{lossless_assumption} does presume that the leakage dissipation of the DC link capacitor is negligible in comparison to the other dissipations in the electrical network.
Relaxation of this assumption to account for a finite leakage resistance would significantly complicate the analysis, because the resultant losses would depend on the DC link voltage.
In the present formulation, this voltage is unnecessary to specify.

\begin{assumption}\label{instantaneous_assumption}
The energy stored in the inductors and capacitor in $\mathcal{N}$ is negligible.
\end{assumption}

This assumption simplifies the analysis, because it implies that the power flowing out of $\mathcal{N}$ may be approximated as zero, for all time.
Let $P_a(t)$ be the power delivered to the plant and to capacitor $C_s$; i.e.,
\begin{equation} \label{eq0_2}
	P_a(t) = \i^T(t) \v(t) + v_s(t) \left( u_s(t)-v_s(t)/R_s \right).
\end{equation}
Then the total power flowing out of $\mathcal{N}$ is $P_a(t)+P_d(t)$. 
The three above assumptions, together, imply that it is zero, i.e., 
\begin{equation} \label{eq0_3} 
	\i^T(t)\v(t) + \i^T(t)\R\i(t) + v_s(t)u_s(t) + R_r u_s^2(t) = 0 , \quad \forall t.
\end{equation}
The reasonableness of Assumption \ref{instantaneous_assumption} depends on the circumstances of the application.
The primary role of the inductors and DC link capacitor is to maintain acceptably-small current and voltage ripple amplitudes over the course of one PWM switching cycle. 
As such, Assumption \ref{instantaneous_assumption} is justifiable if the PWM frequency is sufficiently high, and if the power-electronic hardware components are properly designed.

\begin{assumption}\label{tracking_assumption}
$\mathcal{N}$ facilitates instantaneous tracking between $\i(t)$ and $\i^*(t)$, assuming $\i^*(t)$ is feasible. 
\end{assumption}

Qualitatively, this assumption states that the transient response of the tracking error $\i(t)-\i^*(t)$ is much faster than the time constants of plant $\map{P}$. 
It is justified if the PWM switching frequency is orders of magnitude higher than the bandwidth of $\map{P}$, and enables us to treat the $\i(t)$ as a control input to $\map{P}$.

Analogous models could be adopted for synthetic impedances, rather than admittances.
In this case, the architecture of the SPSA hardware in Fig.~\ref{fig2} has a dual realization, and one may arrive at similar modeling assumptions in this case.
Similarly, the model can be generalized realizations of self-powered immitances (i.e., mixtures of impedance and admittance within the same network).
Additionally, we note that mechanical realizations of SPSA admittances can, in principle, be implemented. 
This could be accomplished using PWM switch-mode hydraulic power trains \cite{van2013fluid_a,van2013fluid_b}. 

\vspace{12pt}

\subsection{Implicit Conditions on SPSA Feasibility}\label{implicit-conditions}

\begin{figure}
	\centering
	\includegraphics[scale=.9]{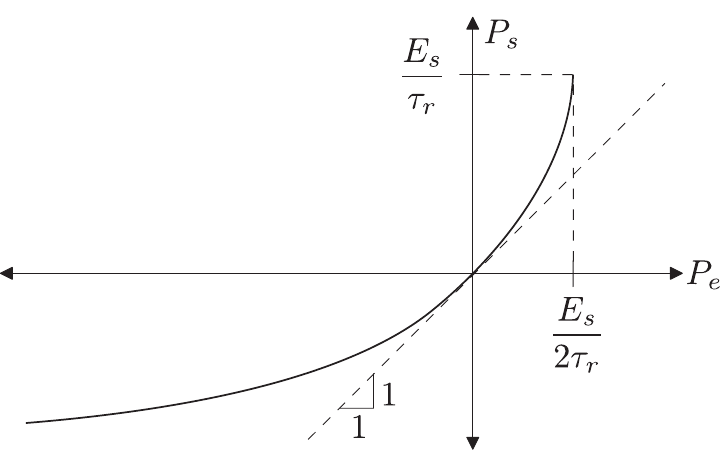}
	\caption{Relationship between actuation power flow ($P_e)$ and power flow to storage ($P_s)$}
	\label{fig4}
\end{figure}

\begin{figure}
	\centering
	\includegraphics[scale=.9]{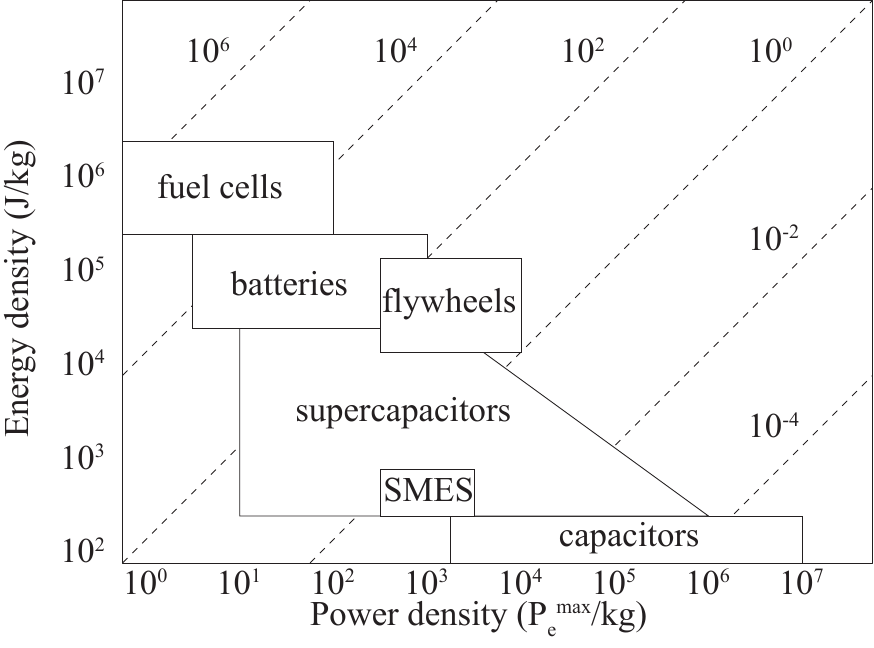}
	\caption{Typical Ragone plot with approximate regions for different types of storage technology: Dotted lines are values of $\tauR$ (in seconds)}
	\label{ragone}
\end{figure}

Quantity $P_s(t) \triangleq -u_s(t)v_s(t)$ is the power evacuated from the energy storage subsystem and sent to the DC link. 
One can view $P_s(t)$ as a ``slack variable'' which takes whatever value is necessary to realize a desired $\i(t)$. 
The power flow quantity 
\begin{equation} \label{eq7}
	P_e(t) = \i^T(t)\v(t) + \i^T(t)\R\i(t)
\end{equation}
is the power network $\mathcal{N}$ must deliver to the plant to realize current $\i(t)$. 
In general $P_e(t) \neq P_s(t)$ because some power (specifically, $P_r(t) = R_r u_s^2(t)$) is dissipated upon transmission from storage to the DC link. 
The energy $E_s$ stored in $C_s$ is 
\begin{equation} \label{ES}
	E_s(t) = \tfrac{1}{2} \C_s v_s^2(t).
\end{equation}
Given $P_e$, constraint \eqref{eq0_3} gives 
\begin{equation}
	u_s(t) = -\tfrac{1}{2R_r} v_s(t) \pm \sqrt{ \left(\tfrac{1}{2R_r}v_s(t)\right)^2 - \tfrac{1}{R_r}P_e(t) }.
\end{equation}
Multiplying by $v_s(t)$, and taking the maximal solution for $u_s(t)v_s(t)$ (i.e., the solution for maximal recharge) gives a relationship between $P_s(t)$, $P_e(t)$, and $E_s(t)$ as
\begin{equation} \label{eq3}
	P_s(t) = \tfrac{1}{\tauR}E_s(t) - 
								\sqrt{ \left(\tfrac{1}{\tauR}E_s(t)\right)^2 - \tfrac{2}{\tauR}E_s(t)P_e(t) }
\end{equation}
where $\tauR \triangleq R_r C_s$.
This relationship is shown in Fig.~\ref{fig4}.  
Note the upper bound $P_e(t) \leqslant E_s/2\tauR$, beyond which $P_s(t)$ has no solution. 
This characterizes the maximum power deliverable from the storage system, as a function of its stored energy, beyond which it will ``choke'' and not be able to supply demanded power.
As such, $\tauR$ characterizes the ability of the storage system to deliver ``pulse power'' to the plant. 
In comparisons of different energy storage technologies, it is common to plot the maximum power density vs. energy density for different storage systems on a log-log scale, which is called a \emph{Ragone plot} \cite{Christen2000}. 
Fig. \ref{ragone} shows an illustration of approximate ranges for different types of storage technologies on this plot. 
Storage systems with the same $\tauR$ lie along the same unit-slope line on the plot.  

Fig.~\ref{fig4} implies that at every time $t$, an SPSA is not at liberty to produce any arbitrary $\i(t)$, even if it has stored up some energy, because there is a limit to the rate at which this stored energy can be delivered to the plant. 
There is therefore an instantaneous feasible domain for $\i(t)$, characterized by
\begin{equation} \label{Pe-const-iv}
	\i^T(t)\R\i(t) + \i^T(t)\v(t) \leqslant \tfrac{1}{2\tauR}E_s(t).
\end{equation}
However, constraint \eqref{Pe-const-iv} is actually a dynamic constraint, because $E_s(t)$ depends on past values of $\i$ and $\v$.
The differential equation for $v_s$ is
\begin{equation} \label{sSdot}
	C_s \tfrac{d}{dt}v_s(t)	= -\tfrac{1}{R_s} v_s(t) + u_s(t).
\end{equation}
Eqs.~\eqref{ES} and \eqref{sSdot} imply
\begin{equation} \label{eq2}
	\tfrac{d}{dt} E_s(t) = -\tfrac{2}{\tauS}E_s(t) + u_s(t)v_s(t)
\end{equation}
where the leakage constant $\tauS \triangleq R_s C_s$.
Eqs.~\eqref{eq2} and \eqref{eq3} imply 
\begin{multline} \label{Econsv2}
	\tfrac{d}{dt} E_s(t) = -\left( \tfrac{2}{\tauS} + \tfrac{1}{\tauR} \right)E_s(t) \\ 
	+ \sqrt{ \left(\tfrac{1}{\tauR}E_s(t)\right)^2 - \tfrac{2}{\tauR} E_s(t) P_e(t) }. 
\end{multline}
In order for control input $\i \in \Leb_{2e}^+$ to be feasible, given $\v \in \Leb_{2e}^+$, it must be such that for any $E_s(0) \in\Reals_{>0}$, differential equation \eqref{Econsv2} has a solution $E_s(t) \in\Reals_{>0}$ for all $t > 0$.

\begin{def1}
Given $\v \in \Leb_{2e}^+$, $\mathbb{U}_{SP}(\v; \R,\tauS,\tauR)$ is the set of all $\i\in\Leb_{2e}^+$, for which \eqref{Econsv2} has a unique solution satisfying $E_s(t)\in\Reals_{>0}$ for all $t>0$, and for each $E_s(0)\in\Reals_{>0}$.
\end{def1}

\begin{def1}
Given $\{\R,\tauR,\tauS\}$, a synthetic admittance $\map{Y}$ as in \eqref{i=-Yv} is called \emph{Self Powered} if $\i\in\mathbb{U}_{SP}(\v;\R,\tauS,\tauR)$, $\forall \v \in \Leb_{2e}^+$.
Given $\{\R,\tauR,\tauS\}$, $\mathbb{Y}_{SP}(\R,\tauS,\tauR)$ is the set of all SPSAs.
\end{def1}


\section{Feasibility of General Linear SPSAs}
\label{sec:sec2}

\subsection{Characterization of linear admittance domain} 

We consider linear SPSAs comprised of a time-varying static gain, together with a Volterra integral.
To state this more precisely, for time-valued matrix $\Y_0(t) \in \Reals^{n_p\times n_p}$ and Volterra kernel $\Y_1(t,\theta) \in \Reals^{n_p\times n_p}$, define the system $\map{S}(\Y_0,\Y_1)$ as
\begin{equation}
(\map{S}\v)(t) \triangleq \Y_0(t) \v(t) + \int_0^t \Y_1(t,\theta) \v(\theta) d\theta.
\end{equation}
In order to ensure certain regularity conditions hold, we make some assumptions on the domains of $\Y_0$ and $\Y_1$.
\begin{def1}
The pair $\{\Y_0,\Y_1\}$ are called \emph{regular} if $\Y_0 \in \Leb_\infty^+$ and $\Y_1 \in \mathbb{D}$ which implies that $\Y_1(t,\cdot) \in \Leb_{2}([0,t]), \forall t > 0$, $\Y_1(\cdot,\theta) \in \Leb_{2e}([\theta,\infty)), \forall \theta > 0$, and $\Y_1 \in \Leb_{2e}(\mathbb{T})$ with 
\begin{equation}
\mathbb{T} = \left\{ (t,\theta) \in \Reals^2 : t \geqslant \theta \geqslant 0 \right\}.
\end{equation}
\end{def1}
Then we can state the domain of linear SPSAs as follows.
\begin{def1} \label{YSPL}
$\mathbb{Y}_{SP}^L(\R,\tauS,\tauR)$ is the largest subset of $\mathbb{Y}_{SP}(\R,\tauS,\tauR)$ for which, for each $\map{Y} \in \mathbb{Y}_{SP}^L(\R,\tauS,\tauR)$, there exists a regular pair $\{\Y_0,\Y_1\}$ such that 
\begin{equation} \label{iv-conv}
\map{Y} = \map{S}(\Y_0,\Y_1).
\end{equation}
\end{def1}

\begin{figure}
\centering
\includegraphics[scale=1]{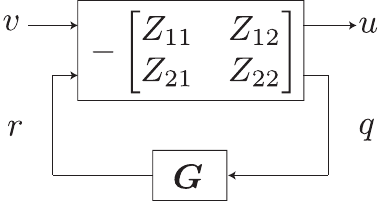}
\caption{LFT representation of linear SPSA $Y$}
\label{lft}
\end{figure}

In this section we derive sufficient conditions on $\{Y_0,Y_1\}$ to guarantee $\map{Y} \in \mathbb{Y}_{SP}^L(\R,\tauS,\tauR)$.
To begin, consider Fig.~\ref{lft}, in which $\map{Y}$ is formulated as a linear fractional transformation (LFT), from a time-varying static gain matrix $\Z(t)$ and a linear time-varying system $\map{G}$ with kernel $\G$, i.e.,
\begin{align} 
\begin{bmatrix} \i(t) \\ \q(t) \end{bmatrix} &= -\Z(t) \begin{bmatrix} \v(t) \\ \p(t) \end{bmatrix} \label{Z-def-ltv} \\
\p(t) &= \int_0^t \G(t,\theta) \q(\theta) d\theta  \label{p-def-ltv}
\end{align}
where $q$ and $r$ are auxiliary signals with $n_r \triangleq \dim(\p) = \dim(\q)$, and where we assume $\{\Z,\G\}$ regular.
We use the shorthand $\map{F}_\ell(\Z,\G)$ to denote the mapping $\v \mapsto -\i$ characterized by \eqref{Z-def-ltv} and \eqref{p-def-ltv}.

The lemma below groups together some useful results related to the LFT in Fig.~\ref{lft}.
For this lemma, we partition $\Z(t)$ as in Fig.~\ref{lft}, where $Z_{11}(t) \in \Reals^{n_p\times n_p}$ and the other partitions are compatibly-sized.

\begin{lm1} \label{LFT-lemma}
\begin{enumerate}
\item
For regular $\{\Z,\G\}$, there exists a unique and regular $\{\Y_0,\Y_1\}$ such that $\map{S}(\Y_0,\Y_1) = \map{F}_\ell(\Z,\G)$. 
\item
Let $\{\Y_0,\Y_1\}$ be regular, and let $\Z\in\Leb_{\infty}^+$ with $\Z_{11}(t) = \Y_0(t)$.
Suppose there exists $\epsilon > 0$ such that $\Z_{12}(t)\Z_{12}(t)^T > \epsilon I$ and $\Z_{21}^T(t)\Z_{21}(t) > \epsilon I$, $\forall t$. 
Then there exists a compatible kernel $\G\in\mathbb{D}$, such that $\map{S}(\Y_0,\Y_1) = \map{F}_\ell(\Z,\G)$.
\item 
If $\{\Z,\G\}$ are regular then $\v\in\Leb_{2e}^+$ $\Rightarrow$ $\i, \q, \p \in \Leb_{2e}^+$
\end{enumerate}
\end{lm1}

\begin{proof}
({\itshape a})
$\Y_0$ and $\Y_1$ are related to $\Z$ and $\G$ via
\begin{align}
\Y_0(t) &= \Z_{11}(t) , & 
\Y_1(t,\theta) &= \Z_{12}(t) T(t,\theta) \label{Y1=-Z12*T}
\end{align}
where resolvent kernel $T(t,\theta)$ is the solution to 
\begin{equation} \label{T-eq}
T(t,\theta) = \G(t,\theta)\Z_{21}(\theta) + \int_\theta^t \G(t,\phi)\Z_{22}(\phi) T(\phi,\theta) d\phi
\end{equation}
over $\{t,\theta\} \in \mathbb{T}$.
Given $\Z$ and $\G$, equation \eqref{T-eq} is a Volterra equation of second kind, in the unknown $T$. 
It is a standard result (see, e.g., \cite{kanwal2013linear}) that a unique solution exists in $\mathbb{D}$ if $\G(t,\phi)\Z_{22}(\phi)$ and $\G(t,\phi)\Z_{21}(\phi)$ are in $\mathbb{D}$.
These conditions are guaranteed if $\G \in \mathbb{D}$, and if $\Z_{22}, \Z_{21} \in \Leb_\infty^+$, both of which are assumed.
Furthermore, if $T\in\mathbb{D}$ then it follows that $\Y_1$ in \eqref{Y1=-Z12*T} is as well, if $\Z_{12}(t)\in\Leb_{\infty}^+$, which is true if $\Z\in\Leb_{\infty}^+$.

({\itshape b})
If $\Z_{12}(t)\Z_{12}^T(t) > \epsilon I$, $\forall t$, then $\Z_{12}(t)$ has a pseudoinverse $\Z^\dagger_{12}(t)$ with $\|\Z^\dagger_{12}(t)\|_2 < \epsilon^{-1/2}$, $\forall t$.
Thus, we conclude that $\Z^\dagger_{12}(t) \in \Leb_{\infty}^+$.
Similarly, if $\Z_{21}(t)\Z_{21}^T(t) > \epsilon I$, $\forall t$, then it has a pseudoinverse $\Z_{21}^\dagger(t) \in \Leb_\infty^+$.
Given $\Y_1\in\mathbb{D}$, \eqref{Y1=-Z12*T} gives a compatible $T$ as
$T(t,\theta) = \Z_{12}^\dagger(t) \Y_1(t,\theta)$.
The facts that $\Y_1 \in \mathbb{D}$ and $\Z_{12}^\dagger \in \Leb_\infty^+$ imply that $T\in\mathbb{D}$.
Hypothesize a solution for kernel $\G(t,\theta) = \Z_{12}^T(t) \bar{\G}(t,\theta) \Z_{21}^\dagger(\theta)$, where we note that $\G\in\mathbb{D}$ if $\bar{\G}$ is.
Such a solution exists if, via \eqref{T-eq}, $\bar{\G}$ is the solution to 
\begin{equation*} 
\bar{\G}(t,\theta) =  \Z_{21}^\dagger(t) T(t,\theta)
- \int_\theta^t \bar{\G}(t,\phi)\Z_{A}(\phi) T(\phi,\theta) d\phi
\end{equation*}
where $Z_A(\phi) \triangleq Z_{21}^\dagger(\phi)Z_{22}(\phi)$.
This is a Volterra integral equation of second kind, and has a unique solution $\bar{\G}\in\mathbb{D}$ if 
$\Z_{21}^\dagger(t)\Z_{22}(t)T(t,\theta)\in\mathbb{D}$ and $\Z_{21}^\dagger(t)T(t,\theta)\in\mathbb{D}$.
By similar arguments as made previously, these conditions are guaranteed by $T\in\mathbb{D}$ and boundedness of $\Z_{21}^\dagger$ and $\Z_{22}$.  

({\itshape c})
Put $\q$ directly in terms of $\v$ as
\begin{equation*}
\q(t) = -\Z_{21}(t) \v(t) - \int_0^t \Z_{22}(t) \G(t,\theta) \q(\theta) d\theta
\end{equation*}
which is a Volterra equation of second kind, and has a unique solution in $\Leb_{2}([0,t_1])$ for any finite $t_1$ if $\Z_{21}(t)\v(t)$ is absolutely integrable over $t\in[0,t_1]$ and if $\Z_{22}(t) \G(t,\theta)$ meets the conditions of regularity for domain $0\leqslant \theta \leqslant t \leqslant t_1$.  
The first of these conditions is guaranteed for any finite $t_1$ as a result of the facts that $\Z_{21}(t)$ is bounded and $\v \in \Leb_{2e}^+$, because any function in $\Leb_2([0,t_1])$ is also in $\Leb_1([0,t_1])$ for finite $t_1$.  
The second condition is guaranteed by the facts that $\Z_{22}(t)$ is bounded and $\G$ is regular.
But if $\q(t)$ over $t\in[0,t_1]$ is in $\Leb_2([0,t_1])$ for each $t_1\geqslant 0$ then we can conclude that $\q(t)$ over $t\in[t_1,t_2]$ is in $\Leb_2([t_1,t_2])$, for arbitrary $0\leqslant t_1 < t_2 < \infty$. 
Thus we conclude that $\q\in\Leb_{2e}^+$.
With this proved, we then note that the solution for $\p$ over interval $t\in[0,t_1]$ for finite $t_1$, as solved by \eqref{p-def-ltv}, must be bounded at each time due to the Cauchy-Schwartz inequality, and if it is bounded, then it is in $\Leb_2([0,t_1])$. 
By a similar reasoning, we conclude that $\p$ over $t\in[t_1,t_2]$ must be in $\Leb_2([t_1,t_2])$ for arbitrary $0\leqslant t_1 <  t_2<\infty$, and therefore that $\p\in\Leb_{2e}^+$.
Finally, we have that $\i(t) = -\left(\Z_{11}(t) \v(t) + \Z_{12}(t) \p(t)\right)$, which is in $\Leb_{2e}^+$ if $\v$ and $\p$ are, because the $\Z_{ij}(t)$ terms are bounded.
\end{proof}

From the above lemma we obtain the following corollary, the proof of which is immediate.
\begin{cor1}
For all $\map{Y} \in \mathbb{Y}_{SP}^L(\R,\tau_s,\tau_r)$ there exist a regular pair $\{\Z,\G\}$ such that $\map{Y} = \map{F}_\ell(\Z,\G)$.
\end{cor1}

\subsection{Derivation of SPSA feasibility conditions}

Below, we derive sufficient conditions on $\{\Z,\G\}$ to ensure that $\map{F}_\ell(\Z,\G) \in \mathbb{Y}_{SP}(\R,\tauS,\tauR)$. 
In presenting the results, the following definition are useful.
\begin{def1} Let $\{t,\theta\}\in\mathbb{T}$, $\alpha\in\Reals$ and $\Gamma\in\mathbb{D}$.
Then $\map{L}_{\Gamma,\alpha}: \Leb_{2e}^+ \rightarrow \Leb_{2e}(\mathbb{T})$ is defined as
\begin{multline}
(\map{L}_{\Gamma,\alpha}\q)(t,\theta) \triangleq  \int_0^\theta \Gamma(\theta,\phi)\q(\phi) d\phi + \int_\theta^t \Gamma^T(\phi,\theta) \q(\phi) d\phi \\
- \alpha\int_0^t \Gamma^T(t,\theta) \Gamma(t,\phi) \q(\phi) d\phi.
\end{multline}
Note that for each $t > 0$, $\map{L}_{\Gamma,\alpha}(t,\cdot) : \Leb_2([0,t]) \rightarrow \Leb_2([0,t])$ is a self-adjoint linear operator.
Consequently, $\spec\left( \map{L}_{\Gamma,\alpha}(t,\cdot) \right) \subseteq \Reals$, and we define 
\begin{equation}
\underline{\lambda} \left( \map{L}_{\Gamma,\alpha}(t,\cdot) \right) \triangleq
\inf \left\{ \spec\left( \map{L}_{\Gamma,\alpha}(t,\cdot) \right) \right\}.
\end{equation}
\end{def1}

\begin{th1} \label{LTV-feas-gen}
$\map{Y} \in \mathbb{Y}_{SP}^L(\R,\tauS,\tauR)$ if there exists a regular pair $\{\Z,\G\}$ such that $\map{Y} = \map{F}_\ell(\Z,\G)$, and such that 
\begin{align}
& \Z(t)+\Z^T(t)-2\Z^T(t)W\Z(t) \geqslant 0 \ , \ \forall t\geqslant 0     \label{Z(t)-norm-condn}
\\
& \underline{\lambda}(\map{L}_{\tilde{G},\tauR}(t,\cdot)) \geqslant 0 \ , \ \forall t > 0		\label{G(t)-def-condn}
\end{align}
where $W \triangleq \text{blockdiag}\{\R,\eye\}$ and $\tilde{G}(t,\theta) \triangleq e^{(t-\theta)/\tauS}G(t,\theta)$.
\end{th1}

The proof of Theorem \ref{LTV-feas-gen} involves three basic steps, which are isolated in the following lemmas.

\begin{lm1} \label{LTV-proof-lm1}
If $\Z$ satisfies \eqref{Z(t)-norm-condn} for $\R\geqslant 0$, then for all $\v \in \Leb_{2e}^+$, the following statement is true:
\begin{equation}
E_s(t) \geqslant \tfrac{\tauR}{2} \|\p(t)\|_2^2 \ , \forall t\geqslant 0  \ \ \Rightarrow \ \ \i \in \mathbb{U}_{SP}(\v;\R,\tauS,\tauR).
\end{equation}
\end{lm1}
\begin{proof}
In order for $\i \in \mathbb{U}_{SP}(\v;\R,\tauS,\tauR)$, \eqref{Pe-const-iv} must hold for all $t \in \Reals_{\geqslant 0}$.
But 
\begin{align*}
& P_e(t) - \tfrac{1}{2\tauR}E_s(t)  
\leqslant
P_e(t) + \|\q(t)+\tfrac{1}{2}\p(t)\|_2^2 - \tfrac{1}{2\tauR}E_s(t) 
\\
&\quad\quad\quad
= -\begin{bmatrix} \v(t) \\ \p(t) \end{bmatrix}^T Q(t)  \begin{bmatrix} \v(t) \\ \p(t) \end{bmatrix} + \tfrac{1}{2\tauR} \left( \tfrac{\tauR}{2}\|\p(t)\|_2^2 - E_s(t) \right)
\end{align*}
where $
Q(t) \triangleq -\Z^T(t) W \Z(t) + \tfrac{1}{2} \Z^T(t) + \tfrac{1}{2} \Z(t) 
$.
We conclude that if $Q(t) \geqslant 0$ (i.e., \eqref{Z(t)-norm-condn} holds) then 
\begin{align*}
P_e(t) - \tfrac{1}{2\tauR}E_s(t) 
\leqslant
\tfrac{1}{2\tauR} \left( \tfrac{\tauR}{2}\|\p(t)\|_2^2 -E_s(t) \right)
\end{align*}
which completes the proof.
\end{proof}

\begin{lm1}\label{LTV-proof-lm2}
Let $V_s(t)$ be defined as
\begin{equation} \label{Vs-def}
V_s(t) \triangleq \int_0^t e^{2(\theta-t)/\tauS} \q^T(\theta) \p(\theta) d\theta .
\end{equation}
Then $\G$ satisfies \eqref{G(t)-def-condn} if and only if for all $\v\in\Leb_{2e}^+$,
\begin{equation} \label{Vs-Es-ineq}
V_s(t) \geqslant \tfrac{\tauR}{2} \|\p(t)\|_2^2 \ \ , \ \ \forall t \geqslant 0.
\end{equation}

\end{lm1}
\begin{proof}
Combining \eqref{Vs-def} and \eqref{Vs-Es-ineq} gives the condition
\begin{equation*}
J \triangleq \int_0^t  e^{2(\theta-t)/\tauS} \q^T(\theta) \p(\theta) d\theta - \tfrac{\tauR}{2} \|\p(t)\|_2^2 \geqslant 0.
\end{equation*}
Substituting \eqref{p-def-ltv} gives 
\begin{multline*} 
J = \int_0^t \int_0^t \tilde{\q}^T(\theta) \big( \tilde{\G}(\theta,\phi) 
- \tfrac{\tauR}{2} \tilde{\G}^T(t,\theta) \tilde{\G}(t,\phi) \big) \tilde{\q}(\phi) d\phi d\theta
\end{multline*}
where $\tilde{q}(\theta) \triangleq e^{(\theta-t)/\tauS}\q(\theta)$, and where we take $\tilde{\G}(\theta,\phi) = 0$ for $\phi>\theta$.  
The above is, furthermore, equivalent to
\begin{equation*}
J = \tfrac{1}{2}\int_0^t \tilde{\q}^T(\theta) (\map{L}_{\tilde{G},\tauR}\tilde{\q})(t,\theta) d\theta.
\end{equation*}
Because $\map{L}_{\tilde{G},\tauR}(t,\cdot)$ is self-adjoint, its spectrum is in $\Reals$, and $J \geqslant 0$ for all $q \in \Leb_{2e}^+$ if and only if  \eqref{G(t)-def-condn} holds.
\end{proof}

\begin{lm1} \label{LTV-proof-lm3}
If $\Z$ satisfies \eqref{Z(t)-norm-condn}, $\G$ satisfies \eqref{G(t)-def-condn}, and $E_s(0) > 0$, then for all $\q\in\Leb_{2e}^+$,
$E_s(t) \in \Reals_{>0}$ and  $E_s(t) > V_s(t) \geqslant 0$, $\forall t \geqslant 0$. 
\end{lm1}
\begin{proof}
If $\Z$ obeys \eqref{Z(t)-norm-condn} then from the proof of Lemma \ref{LTV-proof-lm1}, the term in the square root of \eqref{Econsv2} is bounded by
\begin{equation*} 
\left( \tfrac{1}{\tauR} E_s(t) \right)^2 - \tfrac{2}{\tauR} E_s(t) P_e(t)  \geqslant \tfrac{1}{\tauR^2} E_s(t) \left( E_s(t) - \tfrac{\tauR}{2}\|\p(t)\|_2^2 \right) .
\end{equation*}
Because $\G\in\mathbb{D}$, and $\q\in\Leb_{2e}^+$, it follows that $\lim_{t\rightarrow 0} \p(t) = \p(0) = 0$.
Because $E_s(0) > 0$, it therefore follows that over some finite interval following $t=0$, $E_s(t) > \tfrac{\tauR}{2}\|\p(t)\|_2^2$ and the argument in the square root in \eqref{Econsv2} is positive.
Thus we conclude that there exists an interval $t\in[0,t_0]$ over which $E_s(t)\in\Reals_{\geqslant 0}$.
Furthermore, we can take $t_0$ as the first time at which $E_s(t_0) = \tfrac{\tauR}{2}\|\p(t_0)\|_2^2$, if such a finite time exists.
Otherwise, $t_0 = \infty$.
For the remainder of the proof, we assume $t_0$ is finite, and show that this leads to a contradiction.  

For $t\in[0,t_0]$, $E_s(t) = \tfrac{1}{2} C_s v_s^2(t)$ for some unique positive $v_s(t)$.
Define $\epsilon(t)$ as
\begin{align*}
\epsilon(t)	=& u_s(t)v_s(t)-\q^T(t)\p(t) \\
						=& u_s(t)\sqrt{2E_s(t)/C_s} - \q^T(t)\p(t).
\end{align*}
Power constraint \eqref{eq0_3} is then equivalent to
\begin{equation}  
P_e(t) = -\epsilon(t) - \q^T(t)\p(t) - \tauR \frac{ (\epsilon(t)+\q^T(t)\p(t))^2 }{2E_s(t)} .\label{lem3-ltv-eq1}
\end{equation}
In terms of $\epsilon(t)$, \eqref{Econsv2} becomes
\begin{equation*} 
	\tfrac{d}{dt} E_s(t) = -\tfrac{2}{\tauS}E_s(t) + \epsilon(t) + \q^T(t)\p(t)
\end{equation*}
and consequently
\begin{equation*}
  \tfrac{d}{dt} \left( E_s(t) - V_s(t) \right) = -\tfrac{2}{\tauS} \left( E_s(t) - V_s(t) \right) + \epsilon(t).
\end{equation*}
Remembering that $E_s(0) - V_s(0) = E_s(0) > 0$, if it can be shown that $\epsilon(t)\geqslant 0$ whenever $E_s(t) > V_s(t)$, this is sufficient to conclude that $E_s(t) > V_s(t)$, $\forall t >[0,t_0]$.
Because $V_s(t) \geqslant \tfrac{\tauR}{2}\|\p(t)\|_2^2$, this would then imply that there is no finite $t_0$ at which $E_s(t_0) = \tfrac{\tauR}{2} \|\p(t_0)\|_2^2$, thus proving the lemma.

To show that $E_s(t) > V_s(t) \Rightarrow \epsilon(t) \geqslant 0$,  rewrite \eqref{lem3-ltv-eq1} as
\begin{multline} \label{lm3-ltv-eq2}
0 = P_e(t) + \q^T(t)\p(t) + \tauR\frac{ (\q^T(t)\p(t))^2 }{ 2E_s(t) } \\
	  + \frac{\tauR}{E_s(t)} \left( \frac{\epsilon^2(t)}{2} + \left( \q^T(t)\p(t)+\frac{E_s(t)}{\tauR}  \right)\epsilon(t) \right).
\end{multline}
But if $E_s(t) > V_s(t)$ and \eqref{G(t)-def-condn} holds, then
\begin{equation*}
\frac{ (\q^T(t)\p(t))^2 }{ 2E_s(t) } 
\leqslant \frac{ (\q^T(t)\p(t))^2 }{ 2V_s(t) } 
\leqslant \frac{ \|\q(t)\|_2^2 \|\p(t)\|_2^2 }{ 2V_s(t) } 
\leqslant \frac{ \|\q(t)\|_2^2 }{\tauR}
\end{equation*}
where the last inequality is a consequence of Lemma \ref{LTV-proof-lm2}.
Inserting this result into \eqref{lm3-ltv-eq2} gives the inequality
\begin{multline*}
0 \leqslant P_e(t) + \q^T(t)\p(t) + \q^T(t)\q(t) \\
	  + \frac{\tauR}{E_s(t)} \left( \tfrac{1}{2}\epsilon^2(t) + \left( \q^T(t)\p(t) + \frac{E_s(t)}{\tauR} \right)\epsilon(t) \right).
\end{multline*}
If $\Z(t)$ adheres to \eqref{Z(t)-norm-condn} then
\begin{equation*}
P_e(t) + \q^T(t)\p(t) + \q^T(t)\q(t) \leqslant 0.
\end{equation*}
With this result, and because $E_s(t)>0$ by assumption, 
\begin{equation*}
0 \leqslant \tfrac{1}{2}\epsilon^2(t) + \left( \q^T(t)\p(t)+\tfrac{1}{\tauR}E_s(t) \right)\epsilon(t).
\end{equation*}
This implies that there exists $c(t) \in \Reals_{\geqslant 0}$ such that 
\begin{multline}
\epsilon(t) = -\left( \q^T(t)\p(t)+\tfrac{1}{\tauR}E_s(t) \right) 
\\
\pm \sqrt{ \left( \q^T(t)\p(t)+\tfrac{1}{\tauR}E_s(t) \right)^2 + 2c(t) } \geqslant 0.
\end{multline}
The upper root (corresponding to maximum power flow to storage) is always nonnegative, implying $\epsilon(t)\geqslant 0$.
\end{proof}

With Lemmas \ref{LFT-lemma}--\ref{LTV-proof-lm3}, the proof of Theorem \ref{LTV-feas-gen} is almost immediate.
We have assumed $\{\Z,\G\}$ are regular and $\v \in \Leb_{2e}^+$.
By Lemma \ref{LFT-lemma}, these assumptions imply that $\q\in\Leb_{2e}^+$.
By Lemma \ref{LTV-proof-lm3}, this in turn implies that for any $E_s(0) > 0$, the solution $E_s(t) \in \Reals_{\geqslant 0}$, $\forall t$.

\begin{def1}
$\mathbb{Y}_1(\R,\tauS,\tauR)$ is the set of all $\map{Y}\in\mathbb{Y}_{SP}^L(\R,\tauS,\tauR)$ satisfying the conditions of Theorem \ref{LTV-feas-gen}.
\end{def1}

Theorem \ref{LTV-feas-gen} asserts that its conditions are sufficient to infer feasibility; i.e., $\mathbb{Y}_1(\R,\tauS,\tauR) \subset \mathbb{Y}_{SP}^L(\R,\tauS,\tauR)$
It does not claim that the feasibility conditions are necessary.
However, for finite-dimensional systems its conditions become necessary as well as sufficient for special cases in which $\tauR \rightarrow 0$ and $\tauR \rightarrow \infty$. 
This is proved in the Section IV.

\subsection{SPSA feasibility in the time-invariant case}

In the special case that $\map{Y}$ it time-invariant, it is still useful to characterize $\mathbb{Y}_{SP}(\R,\tauS,\tauR)$ in terms of the feasibility of an associated $\{\Z,\G\}$ pair.
In Theorem \ref{LTV-feas-gen}, all that changes for $\Z$ is that constraint \eqref{Z(t)-norm-condn} becomes a static constraint on a time-invariant matrix $\Z_0$.  
However, for $\G$ a more convenient frequency-domain constraint can be used, which relates feasibility of to an associated positive-real constraint.
This is derived in the theorem below. 

\begin{th1}
Let $\map{Y} = \map{F}_\ell(\Z,\G)$ where $\Z(t)=\Z_0$ and $\G\in\mathbb{D}$ satisfies $\G(t,\theta) = \G(t-\theta,0) \triangleq \G_0(t-\theta)$, $\forall \{t,\theta\}\in\mathbb{T}$.  
Suppose that $\G_0\in\Leb_2^+$, is continuous, and that $\G_0(0)$ is finite. 
Then $\map{Y}\in\mathbb{Y}_1(\R,\tauS,\tauR)$ if and only if $\Z_0$ satisfies \eqref{Z(t)-norm-condn}, and if $\hat{\tilde{G}}_0(s) \triangleq \hat{G}_0(s-\tauS^{-1})$ is analytic for $s \in \mathbb{C}_{>0}$ and satisfies
\begin{align}
& \hat{\tilde{G}}_0^H(j\omega) + \hat{\tilde{G}}_0(j\omega) \geqslant 0  \ \ , \ \ \forall \omega\in\Reals \label{inf-dim-lti-G-condition} 
\\
& \tfrac{1}{2} \tauR \left[ G_0(0) + G_0^T(0) \right] \leqslant I. \label{inf-dim-lti-G0-condition}
\end{align}
\end{th1}

\begin{proof}
(Necessity:)  In the time-invariant case $\map{L}_{\tilde{G},\tauR}$ is
\begin{multline*}
(\map{L}_{\tilde{G},\tauR}q)(t,\theta) 
=
\int_0^\theta \tilde{G}_0(\theta-\phi) q(\phi) d\phi \\
+ \int_\theta^t \tilde{G}_0^T(\phi-\theta) q(\phi) d\phi - \tauR \int_0^t \tilde{G}_0^T(t-\theta)\tilde{G}_0(t-\phi) q(\phi) d\phi
\end{multline*}
where $\tilde{G}_0(t) = e^{t/\tauS}G_0(t)$.
Consequently $\map{L}_{\tilde{G},\tauR}(t,\cdot)$ is positive-semidefinite if and only if, for all $\q\in\Leb_2([0,t])$,
\begin{equation} 
0 \leqslant  
2 \int_0^t q^T(\theta) r(\theta) d\theta 
- \tauR \left\| r(t) \right\|_2^2 
\label{inf-dim-lti-proof-1}
\end{equation}
where 
$r(t) \triangleq \int_0^t \tilde{G}_0(t-\theta) q(\theta) d\theta$.
This implies that the first term on the right-hand side of \eqref{inf-dim-lti-proof-1} be nonnegative for all $\q \in \Leb_2([0,t])$ and all $t > 0$.
It is a classical result \cite{anderson2013network} that this is true if and only if $\hat{\tilde{G}}_0(s)$ is analytic for $s \in \mathbb{C}_{>0}$, and if \eqref{inf-dim-lti-G-condition} holds.
To prove the necessity of \eqref{inf-dim-lti-G0-condition}, take $t\searrow 0$ in \eqref{inf-dim-lti-proof-1}. 
Because $\tilde{G}_0(\cdot)$ is continuous from the right for $t \geqslant 0$, we have that \eqref{inf-dim-lti-proof-1} is guaranteed for infinitesimal $t$ only if
\begin{multline*}
\int_0^t q^T(\theta) \left[ \int_0^\theta \tilde{G}_0(0) q(\phi) d\phi \right] d\theta \\
+ 
\int_0^t \left[ \int_0^\phi \tilde{G}_0(0) q(\theta) d\theta \right]^T q(\phi) d\phi \\
-
\tauR
\int_0^t \int_0^t \q^T(\theta) \tilde{G}_0^T(0)\tilde{G}_0(0)\q(\phi) d\theta d\phi \geqslant 0
\end{multline*}
for all $\q\in\Leb_2([0,t])$. 
Fix $q(\theta) = q_0$ for all $\theta \in [0,t]$ and the above requires
$
\tilde{G}_0(0) + \tilde{G}_0^T(0) - \tauR \tilde{G}_0^T(0)\tilde{G}_0(0) \geqslant 0.
$
Completing the square, and noting $\tilde{G}_0(0) = G_0(0)$, gives \eqref{inf-dim-lti-G0-condition}.

(Sufficiency:) For convenience, define
\begin{equation*}
H(t) \triangleq \left\{ 
\begin{array}{lll} \tilde{G}_0(t)   &:& t \geqslant 0 \\
                   \tilde{G}_0^T(-t) &:& t < 0 \end{array} 
\right. .
\end{equation*}
The assumption that $\hat{\tilde{G}}_0$ is positive-real implies that for any $t > 0$ and $\tau>0$, 
\begin{equation*}
\int_0^{t+\tau} \int_0^{t+\tau} \q^T(\theta) H(\theta-\phi) \q(\phi) d\theta d\phi \geqslant 0 
\end{equation*}
for all $\q\in\Leb_2([0,t+\tau])$. 
Breaking this up into integration spans $[0,t)$ and $[t,t+\tau]$, this is equivalent to
\begin{multline*}
\int_0^t \int_0^t \q^T(\theta) H(\theta-\phi) \q(\phi) d\theta d\phi \\
+
\int_0^t \int_t^{t+\tau} \q^T(\theta) \tilde{G}_0(\theta-\phi) \q(\phi) d\theta d\phi \\
+
\int_t^{t+\tau} \int_0^t \q^T(\theta) \tilde{G}_0^T(\phi-\theta) \q(\phi) d\theta d\phi \\
+
\int_t^{t+\tau} \int_t^{t+\tau} \q^T(\theta) H(\theta-\phi) \q(\phi) d\theta d\phi 
\geqslant 0 .
\end{multline*}
For $\theta \in [t,t+\tau]$ choose $\q(\theta) = \frac{1}{\tau}q_0$, and the above becomes
\begin{multline*}
\int_0^t \int_0^t \q^T(\theta) H(\theta-\phi) \q(\phi) d\theta d\phi \\
+
2 q_0^T 
\left[ \frac{1}{\tau} \int_t^{t+\tau}  \int_0^t \tilde{G}_0(\theta-\phi) \q(\phi) d\phi d\theta \right] \\
+2
q_0^T \left[
\frac{1}{\tau^2}
\int_t^{t+\tau} \int_\phi^{t+\tau} \tilde{G}_0(\theta-\phi) d\theta d\phi \right] q_0 
\geqslant 0 .
\end{multline*}
Because $\tilde{G}_0(\psi)$ is continuous from the right for all $\psi \geqslant 0$, we have that as $\tau\searrow 0$, $G_0(\theta-\phi) \rightarrow G_0(t-\phi)$ for all $\theta \in [t,t+\tau]$ and all $t \geqslant \phi$.
So for $\tau$ infinitesimal,
\begin{multline*}
\int_0^t \int_0^t \q^T(\theta) H(\theta-\phi) \q(\phi) d\theta d\phi 
+
2 q_0^T 
 \int_0^t \tilde{G}_0(t-\phi) \q(\phi) d\phi  \\
+
q_0^T \left[ \tfrac{1}{2} \tilde{G}_0(0) + \tfrac{1}{2}\tilde{G}_0^T(0) \right] q_0 
\geqslant 0 .
\end{multline*}
Choose 
\begin{equation*}
q_0 = -\left[ \tfrac{1}{2} \tilde{G}_0(0) + \tfrac{1}{2}\tilde{G}_0^T(0) \right]^{-1} \int_0^t \tilde{G}_0(t-\phi) \q(\phi) d\phi 
\end{equation*}
and the above becomes
\begin{multline*}
\int_0^t \int_0^t \q^T(\theta) H(\theta-\phi) \q(\phi) d\theta d\phi \\
-
\int_0^t \int_0^t \q^T(\theta) \tilde{G}_0^T(t-\theta) \left[ \tfrac{1}{2} \tilde{G}_0(0) + \tfrac{1}{2}\tilde{G}_0^T(0) \right]^{-1} \\
\times  \tilde{G}_0(t-\phi) \q(\phi) d\theta d\phi 
\geqslant 0 .
\end{multline*}
If \eqref{inf-dim-lti-G0-condition} holds then, because $\tilde{G}_0(0) = G_0(0)$, the above implies
\begin{multline*}
\int_0^t \int_0^t \q^T(\theta) H(\theta-\phi) \q(\phi) d\theta d\phi \\
-
\tauR \int_0^t \int_0^t \q^T(\theta) \tilde{G}_0^T(t-\theta) \tilde{G}_0(t-\phi) \q(\phi) d\theta d\phi \geqslant 0 
\end{multline*}
which is \eqref{inf-dim-lti-proof-1}. 
\end{proof}

\section{Finite-Dimensional Linear SPSA Feasibility}\label{sec:sec3}

We narrow the focus to finite-dimensional state space models for $\map{Y}$, of the form
\begin{align}
\tfrac{d}{dt}\x(t) =& \A(t) \x(t) + \B(t) \v(t)   \label{xdot=Ax+Bv} \\
      -\i(t) =& \C(t) \x(t) + \D(t) \v(t)   \label{i=Cx+Dv}
\end{align}
where $A(t)\in\Reals^{n\times n}$, $B(t)\in\Reals^{n\times n_p}$, $C(t)\in\Reals^{n_p\times n}$, and $D(t)\in\Reals^{n_p\times n_p}$.
An equivalent model is \eqref{Z-def-ltv} and \eqref{p-def-ltv}, with 
\begin{align}
\tfrac{d}{dt} \x(t) =& \A_G(t) \x(t) + \B_G(t) \q(t)  \label{x=AGx+Bq}\\
      \p(t) =& \C_G(t) \x(t)                  \label{p=CGx}
\end{align}
and where
\begin{align}
\A(t) =& \A_G(t) - \B_G(t) \Z_{22}(t) \C_G(t)  \label{AG-BG-CG-Z22-to-A} \\
\B(t) =& -\B_G(t) \Z_{21}(t)                    \label{BG-Z21-to-B} \\
\C(t) =& \Z_{12}(t) \C_G(t)                    \label{CG-Z12-to-C}  \\
\D(t) =& \Z_{11}(t).                            \label{Z11-to-D} 
\end{align}
The associated kernel $\G(t,\theta)$ is $\C_G(t)\Phi_G(t,\theta)\B_G(\theta)$, where $\Phi_G(t,\theta)$ is the state transition matrix of $\A_G(t)$.
Note that \eqref{p=CGx} contains no feed-through term, as this would violate $\G\in\mathbb{D}$.

\subsection{Sufficiency conditions for the time-varying case}

In order to find a finite-dimensional version of Theorem \ref{LTV-feas-gen}, Lemma \ref{LTV-proof-lm2} should be put in terms of $\{\A_G,\B_G,\C_G\}$.  
The true equivalent of Lemma \ref{LTV-proof-lm2}, so derived, results in a more restricted version of the Positive Real Lemma for finite-dimensional time-varying systems. 
More specifically, in place of \eqref{G(t)-def-condn}, the associated feasibility condition is the existence of a transient solution $\P(t)$ to a passivity-type Riccati equation which is bounded at each time $t$ by $\P(t) \geqslant \tauR \C_G^T(t) \C_G(t)$.  
This result follows in a manner which is almost identical to that for the Bounded-Real Lemma for time-varying systems, which has been formulated using Hamilton-Jacobi techniques and proved using conjugate-point arguments \cite{green2012linear}. 

However in this paper, we make use of multiplier techniques to develop more conservative conditions.
The reasoning behind this is merely that this approach leads more immediately to a convenient feasibility condition for $\map{Y} \in \mathbb{Y}_{SP}^L(\R,\tauS,\tauR)$ which can be expressed as a system of matrix inequalities, and because it is equivalent to the aforementioned approach when the system is time-invariant. 

\begin{lm1}\label{ltv-finite-dim-lemma1}
For a mapping $\map{G} : \q \mapsto \p$ as in \eqref{x=AGx+Bq} and \eqref{p=CGx}, the associated kernel $\G$ satisfies \eqref{G(t)-def-condn} if there exists $\P(t)=\P^T(t) > 0$ such that for all $t \in \Reals_{\geqslant 0}$, 
\begin{align}
\tfrac{d}{dt} \P(t) + \tilde{\A}_G^T(t) \P(t) + \P(t) \tilde{\A}_G(t) \leqslant 0  \label{ltv-finite-dim-G-feas1} \\
\C_G(t) = \B_G^T(t) \P(t)                                \label{ltv-finite-dim-G-feas2} \\
\tauR \C_G^T(t)\C_G(t) \leqslant \P(t)                   \label{ltv-finite-dim-G-feas3} 
\end{align}
where $\tilde{\A}_G(t) = \A_G(t) + \tauS^{-1} \eye$.
\end{lm1}
\begin{proof}
We begin by finding conditions for the case with $\tauS = \infty$.
Define $V_s(t)$ as \eqref{Vs-def}. It is then a standard result from robust control (see, e.g., \cite{forbes2010passive, forbes2011linear
}) that for some $\P(t) > 0$, $V_s(t) \geqslant \tfrac{1}{2} \x^T(t) \P(t) \x(t)$ for all $t\in\Reals_{\geqslant 0}$ and all $\q\in\Leb_{2e}^+$ if both \eqref{ltv-finite-dim-G-feas1} and \eqref{ltv-finite-dim-G-feas2} hold for all $t\in\Reals_{\geqslant 0}$.
If these conditions hold together with \eqref{ltv-finite-dim-G-feas3}, then $V_s(t) - \tfrac{\tauR}{2}\|r(t)\|_2^2 \leqslant 0$.
Via Lemma \ref{LTV-proof-lm2}, we conclude that $\G$ satisfies \eqref{G(t)-def-condn}.

Now, to get the case for general $\tauS \in \mathbb{R}_{>0}$, let $\tilde{q}(t) \triangleq e^{t/\tauS} q(t)$, $\tilde{r}(t) \triangleq e^{t/\tauS} r(t)$, and $\tilde{V}_s(t) \triangleq e^{2t/\tauS} V_s(t)$. 
It follows that 
\begin{equation*}
\tilde{V}_s(t) = \int_0^t \tilde{q}^T(\theta) \tilde{r}(\theta) d\theta
\end{equation*}
and $V_s(t) \geqslant \tfrac{\tauR}{2} \|r(t)\|_2^2$ if and only if $\tilde{V}_s(t) \geqslant \tfrac{\tauR}{2} \|\tilde{r}(t)\|_2^2$. 
As such, we conclude that the state space for the mapping $\tilde{q}\mapsto \tilde{r}$ must satisfy the same criterion as derived above for the $\tauS = \infty$ case. 
One realization for this mapping is
\begin{align*}
\tfrac{d}{dt} \tilde{x}(t) =& \tilde{A}_G(t) \tilde{x}(t) + B_G \tilde{q}(t) \\
\tilde{r}(t) =& C_G(t) \tilde{x}(t)
\end{align*}
and therefore by Lemma \ref{LTV-proof-lm2}, conditions \eqref{ltv-finite-dim-G-feas1}, \eqref{ltv-finite-dim-G-feas2}, and \eqref{ltv-finite-dim-G-feas3} imply that $\G$ satisfies \eqref{G(t)-def-condn}.
\end{proof}

\begin{lm1}
Let $\map{Y}=\map{F}_\ell(\Z,\G)\in\mathbb{Y}_1(\R,\tauS,\tauR)$ have a kernel $\G(t,\theta)$ satisfying the conditions of Lemma \ref{ltv-finite-dim-lemma1}, with $\dim(\x) = n$ and state space parameters $\{\A_G(t),\B_G(t),\C_G(t)\}$. 
If $\rank\{\B_G(t)\} < n$ for $t\in\mathbb{J}\subseteq\Reals_{\geqslant 0}$ then if inequality \eqref{ltv-finite-dim-G-feas3} is satisfied strictly, there exists an equivalent pair $\{\check{\Z},\check{\G}\}$ such that $\map{F}_\ell(\check{\Z},\check{\G}) = \map{F}_\ell(\Z,\G)$, and such that the state space for kernel $\check{\G}$ has $\rank\{\check{\B}_G(t)\} = n, \ \forall t\geqslant 0$.
\end{lm1}
\begin{proof}
We begin by showing that if $\B_G(t)$ has $n_q<n$ columns, an equivalent feasible pair $\{\check{\Z},\check{\G}\}$ exists with the same state space dimension, and which has square $\B_G(t)$. 
Specifically, define $\Upsilon = \begin{bmatrix} I_{n_q} & 0_{n_q,n-n_q} \end{bmatrix}$ and then define the parameters for this equivalent system as
\begin{align*}
\check{\Z}(t) =& \tilde{\Upsilon}^T \Z(t) \tilde{\Upsilon}, &
\check{\A}_G(t) =& \A_G(t), \\
\check{\B}_G(t) =& \B_G(t)\Upsilon , &
\check{\C}_G(t) =& \Upsilon^T \C_G(t).
\end{align*}
where $\tilde{\Upsilon} \triangleq \text{blockdiag}\{I,\Upsilon\}$.
It is straightforward to verify that  $\map{F}_\ell(\Z,\G) = \map{F}_\ell(\check{\Z},\check{\G})$, and that constraints \eqref{Z(t)-norm-condn}, \eqref{ltv-finite-dim-G-feas1}, \eqref{ltv-finite-dim-G-feas2}, and \eqref{ltv-finite-dim-G-feas3} are inherited from $\{\Z,\G\}$ to $\{\check{\Z},\check{\G}\}$. 
Thus, for the remainder of the proof we assume $n_q \geqslant n$.

For each $t$ at which $\B_G(t)$ is rank deficient, we have the singular value decomposition
\begin{equation*}
\B_G(t) = \begin{bmatrix} U_1 & U_2 \end{bmatrix} \begin{bmatrix} \Sigma_{11} & 0 \\ 0 & 0 \end{bmatrix} \begin{bmatrix} V_1 \\ V_2 \end{bmatrix}
\end{equation*}
where $U$ is unitary, $\Sigma_{11}>0$ is diagonal, and $V \in \Reals^{n\times n_q}$ satisfies $VV^T = I$.
Let $\delta \in \Reals_{>0}$ be an infinitesimal constant, and define system $\check{\map{G}}$ via the state space parameters $\check{\A}_G(t) = \A_G(t)$, $\check{\B}_G(t) = \B_G(t) + \delta U_2 V_2$, and $\check{\C}_G(t) = \check{\B}_G^T(t) \P(t)$. 
Define $\check{\Z}(t) \triangleq T_1(t) \Z(t) T_1(t)$,
where $T_1(t) = \text{blockdiag}\{I,V_1^T(t)V_1(t)\}$. 
Then it is immediate that $\map{F}_\ell(\Z,\G) = \map{F}_\ell(\check{\Z},\check{\G})$, and that $\check{\G}$ inherits conditions \eqref{ltv-finite-dim-G-feas1} and \eqref{ltv-finite-dim-G-feas2} from $\G$. 
If inequaltiy \eqref{ltv-finite-dim-G-feas3} holds strictly, then it is also inherited for $\check{\G}$ if $\delta > 0$ is sufficiently small.
Thus we conclude that a $\delta > 0$ exists such that the conditions in Lemma \ref{ltv-finite-dim-lemma1} are met for $\check{\G}$.
It remains to show that condition \eqref{Z(t)-norm-condn} for $\check{\Z}(t)$ is inherited from $\Z(t)$, which is shown by
\begin{align*}
&\Z(t)+\Z^T(t)-2\Z^T(t)W\Z(t) \geqslant 0 \nonumber \\ 
&\Rightarrow T_1(t) \left( \Z(t)+\Z^T(t)-2\Z^T(t)W\Z(t) \right) T_1(t) \geqslant 0 \\
&\Rightarrow \check{\Z}(t)+\check{Z}^T(t)-2\check{\Z}^T(t)W\check{\Z}(t) \geqslant 0
\end{align*}
where we used the fact that $W \geqslant T_1(t) W T_1(t)$. 
\end{proof}

With the above lemma, we restrict our attention to realizations $\map{G}$ with $\rank\{\B_G(t)\}=n$, $\forall t\geqslant 0$, with the understanding that the resultant domain also encompasses special cases where $\B_G(t)$ is rank-deficient, and inequality \eqref{ltv-finite-dim-G-feas3} is strict.
This allows us to determine convenient feasibility conditions directly in terms of $\{\A(t),\B(t),\C(t),\D(t)\}$.

\begin{th1}\label{main_theorem}
Let \eqref{xdot=Ax+Bv} and \eqref{i=Cx+Dv} be a realization for $\map{Y}$, with $\dim(\x)=n$.
Then $\map{Y}\in\mathbb{Y}_1(\R,\tauS,\tauR)$ if there exist $\P(t)=\P^T(t) > 0$ and $\X(t)\in\Reals^{n\times n}$  such that
\begin{align}
& \tfrac{d}{dt} P + \A^T \P + \P \A + 2\tauS^{-1}\P + \X + \X^T \leqslant 0   \label{ltv-finite-dim-theorem-lmi1} \\
&\begin{bmatrix} 
-2(\X+\X^T)  & 2\P \B   & 2\C^T         & -2\X^T   \\
 \star       & -\R^{-1} & 2\D^T-\R^{-1} & 2\B^T \P \\
 \star       & \star    & -\R^{-1}      & 0           \\
 \star       & \star    & \star         & -\tfrac{1}{\tauR}\P  
\end{bmatrix} 
\leqslant 0 
\label{ltv-finite-dim-theorem-lmi2}
\end{align}
hold for each $t\in\Reals_{\geqslant 0}$, where in the above equations we have suppressed the dependency of $A, B, C, D, X,$ and $P$ on $t$. 
\end{th1}
\begin{proof}
In the interest of brevity, we suppress explicit notation of time-dependence throughout.
Define $\X$ as
\begin{equation}\label{ltv-finite-dim-theorem-proof-2}
\X = \P \B_G \Z_{22} \B_G^T \P 
\end{equation}
and \eqref{ltv-finite-dim-G-feas1} becomes \eqref{ltv-finite-dim-theorem-lmi1}, using \eqref{AG-BG-CG-Z22-to-A} and \eqref{ltv-finite-dim-G-feas2}. 
To get \eqref{ltv-finite-dim-theorem-lmi2}, first note that for $\R > 0$ (and therefore $W>0$), \eqref{Z(t)-norm-condn} is equivalent, through a Schur complement, to
\begin{equation} \label{ltv-finite-dim-theorem-proof-1}
\begin{bmatrix} -W^{-1} & 2\Z^T - W^{-1} \\ 2\Z - W^{-1} & -W^{-1} \end{bmatrix} \leqslant 0 .
\end{equation}
Define $T =\blkdiag\{ I, \P\B_G , I , \P\B_G \}$ and note that because we assume $\rank\{\B_G\} = n$ and because $\P > 0$, it follows that $T^T T > 0$. As such, \eqref{ltv-finite-dim-theorem-proof-1} is true if and only if
\begin{equation*}
T \begin{bmatrix} -W^{-1} & 2\Z^T - W^{-1} \\ 2\Z - W^{-1} & -W^{-1} \end{bmatrix} T^T \leqslant 0.
\end{equation*}
Multiplication, and use of \eqref{BG-Z21-to-B}, \eqref{CG-Z12-to-C}, \eqref{Z11-to-D}, \eqref{ltv-finite-dim-G-feas2}, and \eqref{ltv-finite-dim-theorem-proof-2} gives
\begin{equation*}
\begin{bmatrix} Q_1 & Q_2^T \\ Q_2 & Q_1 \end{bmatrix} \leqslant 0 
\end{equation*}
where $Q_1 \triangleq \text{blockdiag}\{ -\R^{-1}, -\Xi \}$ and
\begin{align*}
Q_2(t) \triangleq& 
\begin{bmatrix}
2\D-\R^{-1}  & 2\C   \\
-2\P\B        & 2\X - \Xi
\end{bmatrix} 
\end{align*}
with $\Xi \triangleq \P\B_G\B_G^T\P$.
Rearranging rows and columns gives
\begin{equation*}
\begin{bmatrix}
-\Xi          & 0         & 2\C^T         & 2\X^T-\Xi   \\
 \star        & -\R^{-1}  & 2\D^T-\R^{-1} & -2\B^T\P         \\
 \star        & \star          & -\R^{-1}      & 0                  \\
 \star        & \star          &  \star             & -\Xi
\end{bmatrix}
\leqslant 0.
\end{equation*}
We note that $\Xi$ is not implied by a given $\{\A,\B,\C,\D,\P,\X\}$; i.e., it is an independent parameter.
From our assumption that $\rank\{\B_G\} = n$, we can infer that $\Xi$ is invertible. 
Via a Schur complement, the above is equivalent to 
\begin{multline*}
\begin{bmatrix}
-2\X^T-2\X   & 2\P\B    & 2\C^T         \\
 \star       & -\R^{-1} & 2\D^T-\R^{-1} \\
 \star       & \star    & -\R^{-1}      
\end{bmatrix}
\\
+
\begin{bmatrix} -2\X^T \\ 2\B^T\P \\ 0 \end{bmatrix} \Xi^{-1} \begin{bmatrix} -2\X^T \\ 2\B^T\P \\ 0 \end{bmatrix}^T 
\leqslant 0 .
\end{multline*}
Enforcing \eqref{ltv-finite-dim-G-feas3} requires $\Xi \leqslant \tfrac{1}{\tauR}P$. 
The above inequality is made least conservative by setting $\Xi$ equal to this upper bound.
Note that because $P > 0$, this confirms that $\rank\{\B_G\} = n$.
Another Schur complement then gives \eqref{ltv-finite-dim-theorem-lmi2}.
\end{proof}

\subsection{Sufficiency conditions in the time-invariant case}

In the case where $\map{Y} \in \mathbb{Y}_{SP}^L(R,\tauS,\tauR)$ is time-invariant, there exist constant matrices $\{A,B,C,D\}$ which realize it via equations \eqref{xdot=Ax+Bv} and \eqref{i=Cx+Dv}. 
In this case we have the following corollary to Theorem \ref{main_theorem}, the proof of which is immediate.

\begin{cor1}\label{lti_main_theorem}
Let \eqref{xdot=Ax+Bv} and \eqref{i=Cx+Dv} be a realization for $\map{Y}$, with $\dim(\x)=n$, and with $A$, $B$, $C$, and $D$ constant.
Then $\map{Y}\in\mathbb{Y}_1(\R,\tauS,\tauR)$ if there exist time-invariant matrices $\P=\P^T > 0$ and $\X\in\Reals^{n\times n}$ such that \eqref{ltv-finite-dim-theorem-lmi1} and \eqref{ltv-finite-dim-theorem-lmi2} hold. 
\end{cor1}

\subsection{Necessary conditions in the time-invariant case}

Corollary \ref{lti_main_theorem} constitutes a generalization of the Positive Real Lemma. 
Indeed, in the time-invariant case, feasibility conditions \eqref{ltv-finite-dim-theorem-lmi1} and \eqref{ltv-finite-dim-theorem-lmi2} distill to the standard matrix inequality for the Positive Real Lemma, i.e., 
\begin{equation}
\begin{bmatrix}
A^T P + P A & P B - C^T \\ B^TP - C & -D^T-D 
\end{bmatrix} \leqslant 0
\end{equation}
in the case where $\tauR \rightarrow 0$, $\R \rightarrow 0$, and $\tauS \rightarrow \infty$. 
However, the Positive Real Lemma states that the above condition, which must be satisfied for $P=P^T > 0$, is both necessary and sufficient for mapping $\map{Y}$ to be passive. 
This motivates us to investigate more general circumstances in which conditions \eqref{ltv-finite-dim-theorem-lmi1} and \eqref{ltv-finite-dim-theorem-lmi2} are both necessary and sufficient for $\map{Y} \in \mathbb{Y}_{SP}^L(\R,\tauS,\tauR)$. 
In Lemma \ref{necessary_conditions_lemma} below, we derive a necessary condition for $\map{Y} \in \mathbb{Y}_{SP}^L(\R,\tauS,\tauR)$ and then in Theorem \ref{necessary_conditions_theorem}, we show that this condition is equivalent to the sufficient conditions in two special cases.

\begin{lm1}\label{necessary_conditions_lemma}
Let \eqref{xdot=Ax+Bv} and \eqref{i=Cx+Dv} be a realization for $\map{Y}$, with $\dim(\x)=n$, and with $A$, $B$, $C$, and $D$ constant.
Suppose that $\map{Y}\in\mathbb{Y}_{SP}^L(\R,\tauS,\tauR)$ for $\R > 0$, $\tauR \geqslant 0$ and $\tauS > 0$.
Then this necessarily implies the following equivalent conditions:
\begin{enumerate}
\item $\|2\R^{1/2}\hat{\tilde{\Y}}\R^{1/2}-\eye\|_{\mathbb{H}_\infty} \leqslant 1$,
where $\hat{\tilde{\Y}}(s) \triangleq \hat{\Y}(s-\tauS^{-1})$.
\item
There exists $\P=\P^T > 0$ such that 
\begin{equation} \label{necessary_conditions_lemma_LMI}
\begin{bmatrix} 
 2\A^T\P + 2\P\A + 4\tauS^{-1}\P  
             & 2\P\B     & 2\C^T         \\
 \star      & -\R^{-1}  & 2\D^T-\R^{-1} \\
 \star      & \star          & -\R^{-1} 
\end{bmatrix} \leqslant 0 .
\end{equation}
\end{enumerate}
\end{lm1}
\begin{proof}
Let $\v\in\Leb_{2e}^+$ and let $\i = -\map{Y}\v$.
From \eqref{eq3} and \eqref{eq2}, we have that 
\begin{align*}
\tfrac{d}{dt} E_s(t) 
\leqslant & -\tfrac{2}{\tauS} E_s(t) - \i^T(t) \v(t) - \i^T(t) \R \i(t).
\end{align*}
For $\v, \i \in \Leb_{2e}^+$, let $\bar{E}_s(t)$ be the solution to
\begin{align*}
\tfrac{d}{dt} \bar{E}_s(t) 
=& -\tfrac{2}{\tauS} \bar{E}_s(t) - \i^T(t) \v(t) - \i^T(t) \R \i(t)
\end{align*}
with initial condition $\bar{E}_s(0) = E_s(0) > 0$.  It then follows from the Comparison Lemma \cite{khalil2001nonlinear} that $\bar{E}_s(t) \geqslant E_s(t)$, $\forall t > 0$. 
Now, let $\tilde{E}_s(t) \triangleq \bar{E}_s(t) e^{2t/\tauS}$, $\tilde{u}(t) \triangleq e^{t/\tauS} \R^{1/2} \i(t)$, and $\tilde{v}(t) \triangleq e^{t/\tauS} \R^{-1/2} \v(t)$, and note that these quantities are in $\Leb_{2e}^+$ if $\bar{E}_s$, $\tilde{u}$, and $\tilde{v}$ are.
It follows that 
\begin{equation}\label{necessary_proof_1}
\tilde{E}_s(t) = \int_0^t \left( -\tilde{u}^T(\theta) \tilde{v}(\theta) - \tilde{u}^T(\theta) \tilde{u}(\theta) \right) d\theta
\end{equation}
where we note that $\bar{E}_s(t) \geqslant 0, \ \forall t > 0$ if and only if $\tilde{E}_s(t) \geqslant 0, \ \forall t>0$. 
Rearranging \eqref{necessary_proof_1}, the latter of these conditions requires that the transfer function for the mapping $\tilde{v} \mapsto 2\tilde{u}+\tilde{v}$ have an infinity norm no greater than $1$. But it is straightforward to show that if the transfer function for $\v\mapsto \i$ is $-\hat{Y}(s)$, then the transfer function for $\tilde{v} \mapsto 2\tilde{u}+\tilde{v}$ is $-2R^{1/2}\hat{Y}(s-\tauS^{-1})R^{1/2} + I$.
We conclude that condition (a) given in the lemma is necessary for $\map{Y} \in \mathbb{Y}_{SP}^L(\R,\tauS,\tauR)$.  

To prove condition (b), we note that if \eqref{xdot=Ax+Bv} and \eqref{i=Cx+Dv} are a realization for $\map{Y} : \v \mapsto -\i$, with $A$, $B$, $C$, and $D$ constant, then the mapping $\tilde{v} \mapsto -2\tilde{u}+\tilde{v}$ has the realization
\begin{align*}
\tfrac{d}{dt} \x(t) =& \left[ A + \tauS^{-1} I \right] \x(t) + B \R^{-1/2} \tilde{v}(t) \\
\tilde{u}(t) =& -2\R^{1/2} C \x(t) - 2\R^{1/2} D \R^{1/2} + I \tilde{v}(t)
\end{align*}
Via the Bounded Real Lemma \cite{anderson2013network}, condition (a) is true if and only if there exists $P=P^T>0$ such that
\begin{equation*}
\begin{bmatrix}
A^T P + PA + \tfrac{2}{\tauS} P & PBR^{-1/2} & 2C^TR^{1/2} \\
\star & -I & 2R^{1/2} D^T R^{1/2} - I \\
\star & \star & -I
\end{bmatrix}
\leqslant 0
\end{equation*}
which is equivalent to condition (b).
\end{proof}

Matrix inequality \eqref{necessary_conditions_lemma_LMI} in condition (b) of Lemma \ref{necessary_conditions_lemma} is a necessary condition for time-invariant $\map{Y}\in\mathbb{Y}_{SP}^L(\R,\tauS,\tauR)$.
It resembles the two matrix inequalities \eqref{ltv-finite-dim-theorem-lmi1} and \eqref{ltv-finite-dim-theorem-lmi2}, which together constitute sufficient conditions. 
Indeed, as Theorem \ref{necessary_conditions_theorem} below proves, the necessary and sufficient conditions become the same, in the special cases where $\tauR = 0$ and $\tauR \rightarrow \infty$.

\begin{th1}\label{necessary_conditions_theorem}
Let \eqref{xdot=Ax+Bv} and \eqref{i=Cx+Dv} be a minimal realization for $\map{Y}$, with $\dim(\x)=n$, and with $A$, $B$, $C$, and $D$ constant.
Let $\R > 0$, $\tauS > 0$, and either $\tauR = 0$ or $\tauR \rightarrow\infty$. 
Then $\map{Y}\in\mathbb{Y}_{SP}^L(\R,\tauS,\tauR)$ if and only if there exist time-invariant matrices $\P=\P^T > 0$ and $\X\in\Reals^{n\times n}$ such that \eqref{ltv-finite-dim-theorem-lmi1} and \eqref{ltv-finite-dim-theorem-lmi2} hold. 
\end{th1}
\begin{proof}
With $\tauR = 0$, sufficient condition \eqref{ltv-finite-dim-theorem-lmi2} becomes 
\begin{equation*}
\begin{bmatrix} 
-2(\X+\X^T)  & 2\P \B   & 2\C^T         \\
\star        & -\R^{-1} & 2\D^T-\R^{-1} \\
\star       &  \star        & -\R^{-1}      
\end{bmatrix} 
\leqslant 0 
\end{equation*}
for any $P=P^T > 0$ and $X$ finite.
Choosing the least-conservative $\X$ satisfying both this inequality and \eqref{ltv-finite-dim-theorem-lmi1} gives necessary condition \eqref{necessary_conditions_lemma_LMI}. 
With $\tauR \rightarrow \infty$, constraint \eqref{Pe-const-iv} implies that for $E_s(t)$ finite, it is necessary that for all $\v\in\Leb_{2e}^+$ and all $t \geqslant 0$, $P_e(t) \leqslant 0$.
For this to hold for all $\v\in\Leb_{2e}^+$, it must be the case that $\map{Y}$ is a static mapping, i.e., $\hat{Y}(s) = D$, which must satisfy
\begin{equation}\label{necessary_conditions_theorem_proof_1}
\begin{bmatrix} -\R^{-1} & 2D^T - \R^{-1} \\ 2D - \R^{-1} & -\R^{-1} \end{bmatrix} \leqslant 0 .
\end{equation}
Because it is assumed that the realization of $\map{Y}$ in \eqref{xdot=Ax+Bv} and \eqref{i=Cx+Dv} is minimal, it follows that $\dim\{\x\} = 0$, which makes sufficient condition \eqref{ltv-finite-dim-theorem-lmi1} superfluous, and reduces condition \eqref{ltv-finite-dim-theorem-lmi2} to \eqref{necessary_conditions_theorem_proof_1}, implying that it is both necessary and sufficient for $\map{Y} \in \mathbb{Y}_{SP}^L(\R,\tauS,\infty)$. 
\end{proof}

\section{Examples}\label{sec:sec4}

Suppose a linear, passive, colocated controller, characterized by its effective admittance $\map{Y}$, has been designed.
Given this design, one can characterize the domain of parasitic loss parameters for which $\map{Y}$ can be realized with self-powered control.
In other words, we seek to characterize the domain of triples $(R,\tau_s,\tau_r)$ for which a given $\map{Y}$ is a feasible SPSA.
This can be found by considering that for given values of $\tau_s$ and $R$, the value of $\tau_r$ can be maximized subject to the constraint  $\bm{Y} \in \mathbb{Y}_1(R,\tau_s,\tau_r)$.  
Repeating this optimization for all values of $R$ and $\tau_s$ generates a surface bounding the region of feasible parameters in the $(R,\tau_s,\tau_r)$ space. 
This surface is Pareto-optimal, in the sense that it signifies the parameter triples of least efficiency, for which $\map{Y}\in\mathbb{Y}_1(R,\tau_s,\tau_r)$ is assured.
In this section, we present three numerical examples that demonstrate solution techniques for generating these Pareto surfaces.

We emphasize that the primary objective of this section is to characterize feasible loss parameters for a fixed $\map{Y}$.
This is distinct from the much more challenging problem of determining $\map{Y}$, for fixed $(R,\tau_s,\tau_r)$, to optimize some measure of the closed-loop mapping $w \mapsto z$ in Fig. \ref{block}.
The solution to such a constrained optimal control problem appears to be nonconvex for any standard measure of closed-loop performance (such as $\mathbb{H}_2$ or $\mathbb{H}_\infty$ measures) and it is beyond the scope of this paper to provide a numerically-tractable approach to solve such nonconvex optimizations.
Rather, in each example the design for $\map{Y}$ is obtained through a sub-optimal design procedure, which in each case is adapted from a procedure published elsewhere in the open literature.

\subsection{Finite-dimensional, LTI SPSAs}

In this example we consider a self-powered control system to suppress the vibratory response of a civil structure subjected to an earthquake disturbance. Full modeling details can be found in \cite{ligeikis2021nonlinear}, and here we merely summarize its main features.
The plant is a base-isolated, five-story building with a tuned vibration absorber (TVA) located on the top floor. 
The ground acceleration is modeled as stationary filtered white noise with known spectrum \cite{lin1987evolutionary}. 
Control actuation is facilitated using two permanent-magnet synchronous machines, which are mechanically coupled to the rectilinear lateral motion of the structure via ballscrews. 
One actuator imposes a force between the base and the ground, while the second imposes a force between the top story and the TVA. 
They are presumed to have linear back-EMF constants of $1800$V-s/m and $900$V-s/m, respectively, and identical resistances $R_1=R_2=1\Omega$. The energy storage subsystem is assumed to have time constants $\tau_s = 10$s and $\tau_r=0.075$s. 

The performance objective for control design is to minimize the multi-objective LQG performance measure, i.e., 
$
J \triangleq \max_{j\in \{1,...n_z\}} \Ex \{z_j^2\}
$
where $z_j$ are performance outputs of the augmented system, $n_z$ is the number of such outputs, and $\Ex\{\cdot\}$ denotes the expectation in stationarity. Components of $z$ correspond to weighted base displacement, base acceleration, and top-floor acceleration. 
Determination of the $\map{Y} \in \mathbb{Y}_1(\R,\tauS,\tauR)$ which minimizes $J$ constitutes a nonconvex optimization, and there is no known technique to convert it to a convex problem without the introduction of conservatism. 
Here, we make use of the suboptimal, but convex, synthesis procedure described in \cite{ligeikis2021feasibility}.
The procedure consists of two stages. 
In the first stage, a finite-dimensional, strictly-proper, LTI output-feedback controller with state-space parameters $\{A_K,B_K,C_K\}$ is optimized, subject to the relaxed constraint $\Ex\{ P_e(t) \} < 0$.
This constraint is necessary but not sufficient to guarantee self-powered feasibility. 
In the second stage, a feasible SPSA $\map{Y}$ is obtained by setting $A = A_K$ and $B=B_K$, and finding $\{C,D\}$ which optimally project the controller found in the first stage onto $\mathbb{Y}_{1}(\R,\tauS,\tauR)$.

This procedure yields an SPSA $\map{Y}$ with $A$ equal to
\begin{equation}
    A = \text{blockdiag} \left\{ \begin{bmatrix}
    \alpha_i & \beta_i \\ -\beta_i & \alpha_i
    \end{bmatrix}_{i=1,...,3},~\begin{bmatrix}
    \sigma_1 & 0 \\ 0 & \sigma_2 
    \end{bmatrix} \right\}
\end{equation}
where $\alpha_1=-1.676,~\alpha_2=-4.172,~\alpha_3=-7.869,~\beta_1=40.639,~\beta_2=12.225,~\beta_3=3.972,~\sigma_1=-0.287,$ and $\sigma_2=-3.371$. Matrices $B$, $C$, and $D$ are:
\begin{equation}
    B \! =  \begin{bmatrix}
    -0.0661 & -4.6\times10^{-4}\\ -0.0225 & -1.6\times10^{-4}\\ 0.0545 & 4.3\times10^{-4}\\ 0.15 & 9.9\times10^{-4}\\ 0.0222 & 3.4\times10^{-5}\\ -0.394 & -0.00268\\ 0.0244 & 1.8\times10^{-4}\\ -0.223 & -0.00163
    \end{bmatrix}\!,  C^T \!=  \begin{bmatrix}
    0.0749 & 0.0427\\ -0.336 & 0.048\\ -0.0169 & 0.223\\ -0.163 & -0.119\\ -0.309 & 1.02\\ -0.506 & 0.423\\ 0.079 & 0.0782\\ 0.44 & -0.439
    \end{bmatrix}
\end{equation}
\begin{equation}
    D = \begin{bmatrix}
    0.00996 & 0.00237\\ -0.00777 & 0.00464 
    \end{bmatrix}.
\end{equation}

As mentioned, the resultant $\map{Y}$ is sub-optimal, and it is impossible to know the sacrifice in performance (relative to optimality) without solving the original nonconvex problem.
However, it is bounded by the difference in performance between the stage-1 design (which provides a lower bound for the optimal $J$), and that of the stage-2 design.
This is examined in more detail in \cite{ligeikis2021feasibility}.

\begin{figure}
  \centering
	\includegraphics[scale=0.55]{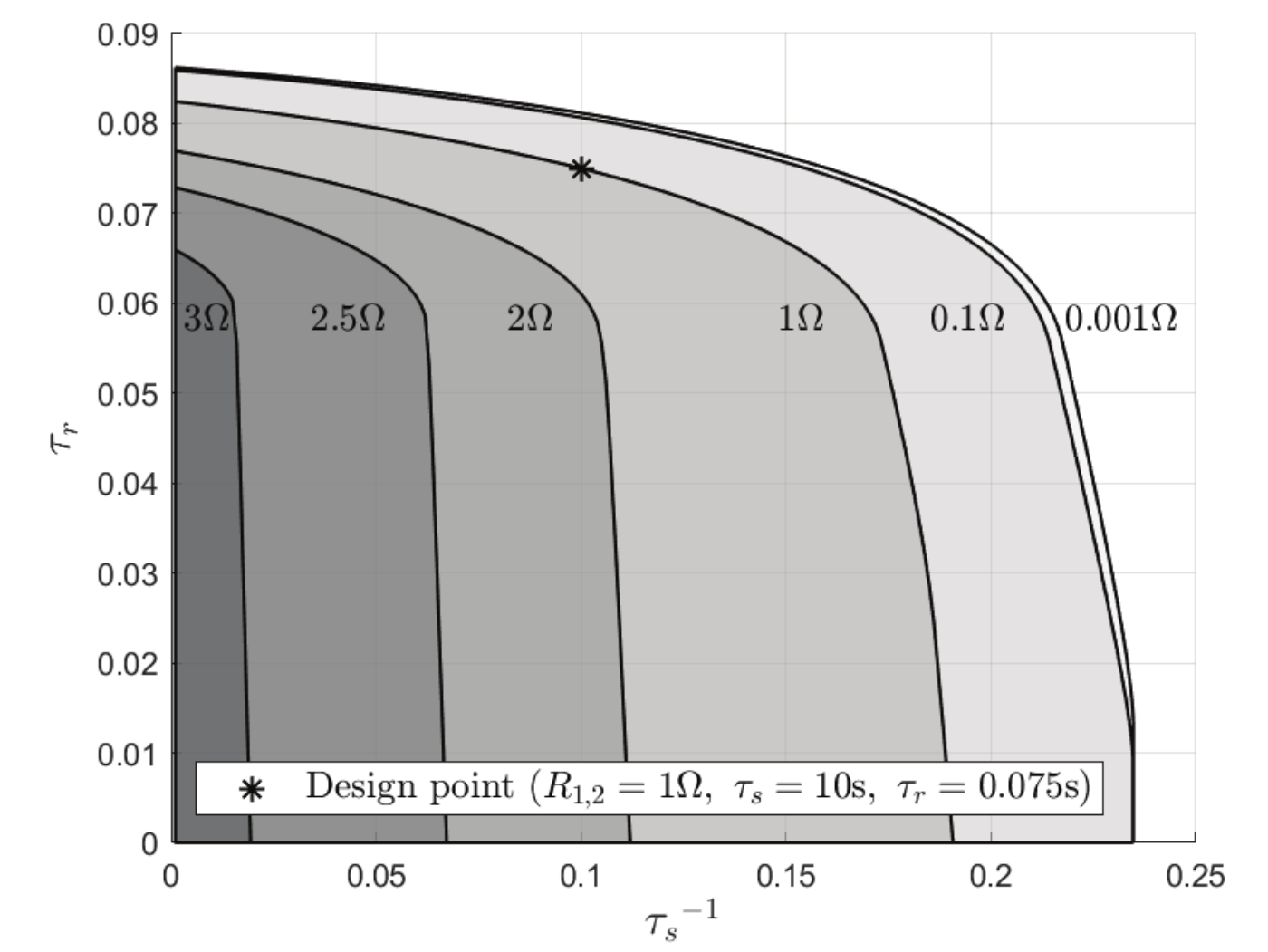}
	\caption{$\{\tau_s^{-1},\tau_r\}$ Pareto fronts (i.e., feasible regions) associated with LTI, SPSA controller, for various values of $\R$ }
	\label{LTI_pareto_fronts}
\end{figure}

We now consider the generation of the Pareto surface of $(R,\tau_s,\tau_r)$ parameters, as described at the beginning of this section, for $\map{Y}$ as given above.
For each $(R,\tau_s)$ pair, the least-efficient (i.e., largest) $\tau_r$ value that renders $\map{Y}$ feasible is found by solving the following optimization problem, based on Corollary \ref{lti_main_theorem}:
\begin{equation*} 
\begin{array}{lll}
\text{Given} &:& A,B,C,D,R,\tau_s \\
\text{Maximize} &:& \tau_r \\
\text{Over} &:& \tau_r,P,X \\
\text{Subject to} &:& \eqref{ltv-finite-dim-theorem-lmi1},\eqref{ltv-finite-dim-theorem-lmi2}, P=P^T>0
\end{array}
\end{equation*}
In this optimization, matrices $P$ and $X$ are restricted to be time-invariant. 
This optimization is quasiconvex, due to bilinear term $\tau_r^{-1}P$, and can be solved efficiently via bisection. 
Fig. \ref{LTI_pareto_fronts} shows the resultant regions of feasible pairs $\{\tau_s^{-1},\tau_r\}$, for various values of $R_1=R_2$. 
The shape of the feasible regions (denoted by the shaded areas in Fig. \ref{LTI_pareto_fronts}) is highly dependent on the particular admittance under study. However, there are some general trends. Clearly, as $R$ decreases, the feasible domain of $\{\tau_s^{-1},\tau_r\}$ pairs increases. In physical terms, this means that as the efficiency of the actuators improves, the energy storage subsystem can be made less efficient, while still realizing the given control law. However, this is only true up to a point, as suggested by the marginal difference in feasible regions for $R_{1,2}=0.1\Omega$ and $R_{1,2}=0.001\Omega$. In addition, as the transmission of energy from storage becomes less efficient (i.e., $\tau_r \rightarrow 0$), the storage system must become more efficient at retaining energy (i.e., $\tau_s \rightarrow \infty)$, and vice versa. 

\subsection{Finite-dimensional, LTV SPSAs}

Next, we synthesize a linear time-varying (LTV), passive controller using a finite-horizon, optimal control framework, and demonstrate how its self-powered feasibility may be assessed. This example is motivated by the plant described in \cite{forbes2010passive}, comprised of a single-degree-of-freedom mechanical system with time-invariant stiffness $k$, time-invariant viscous damping $c$, and mass that decreases in time as $ m(t) = m_f e^{-\alpha t} + m_i$, where $\alpha > 0$ is a constant.  We assume the system is controlled via an electromechanical force transducer, with force $f(t) = K_e u(t)$.
The corresponding internal transducer voltage $v(t) = K_e y(t)$, where $y(t)$ is the mass velocity and $K_e$ is the back-EMF constant. The system dynamics are characterized by the LTV state space
\begin{align}
    \tfrac{d}{dt}x_p(t) =& A_p(t)x_p(t) + B_{pw}(t) w(t) + B_{pu}(t)u(t) \\
    v(t) =& C_p(t)x_p(t)
\end{align}
where $w$ is zero-mean white process noise with unit-intensity. The initial state $x_p(t_i)$ is uncertain with mean $\Ex\{x_p(t_i)\}=0$ and covariance $\Ex\{x_p(t_i)x_p^T(t_i)\}=\Theta$. 
The plant in this example can be shown to be output strictly passive.

We assume $v(t)$ is measured for feedback.
Our objective is to design an passive admittance $\map{Y}$ as in \eqref{i=-Yv}, which is realized via the LTV state space equations \eqref{xdot=Ax+Bv} and \eqref{i=Cx+Dv}.
Then, as discussed at the beginning of this section, we will determine the domain of $(R,\tau_s,\tau_r)$ values for which this $\map{Y}$ is an SPSA.

The methodology used to design $\map{Y}$ was inspired by \cite{forbes2010passive}, but has novel aspects.
We therefore describe it in detail.
To ease the notation, time-dependence of variables are suppressed inside integrals. 
The control objective is to minimize $J$, defined as
\begin{equation}
    J \triangleq \Ex\bigg\{ \int_{t_i}^{t_f} \left(x_p^T Q x_p + u^T N u\right) dt 
    + x_p^T (t_f) S_f x_p(t_f)\bigg\}
\end{equation}
where $Q(t) \geq 0$, $S_f \geq 0$, and $N(t) > 0$, $\forall t \in [t_i,t_f]$.
It is a classical result that $J$ can be expressed as
\begin{equation} 
    J = J_0 + \Ex\left\{\int_{t_i}^{t_f} \left( K x_p - u \right)^T N \left( K x_p - u \right) dt\right\} \label{ltv_perf_obj}
\end{equation}
where 
\begin{equation}
J_0 \triangleq \Ex\left\{ x^T(t_i)S(t_i)x(t_i)\right\} + \textrm{tr}\left\{\int_{t_i}^{t_f} B_{pw}^TSB_{pw}dt \right\},
\end{equation}
$S(t)$ is the solution to the Riccati differential equation
\begin{multline}
    \tfrac{d}{dt}S(t) + A_p^T(t) S(t) + S(t)A_p(t) + Q(t) \\ - S(t)B_{pu}(t)N^{-1}(t)B_{pu}^T(t)S(t) = 0 \label{fsb_riccati}
\end{multline}
with final-value condition $S(t_f)=S_f$, and $K(t)\triangleq -N^{-1}(t) B_{pu}^T(t) S(t)$. 
Note that only the last term in \eqref{ltv_perf_obj} is affected by control input $u$.

We implement a time-varying Luenberger obeserver to estimate $x_p(t)$, and impose the following feedback law
\begin{align}
    \tfrac{d}{dt}\hat{x}_p(t) =& \left[A_p(t) + L(t)C_p(t)\right]\hat{x}_p(t) +B_{pu}(t)u(t) - L(t)v(t) \\
    -u(t) =& \left[F(t)C_p(t)-K(t)\right]\hat{x}_p(t) - F(t) v(t)
\end{align}
where $L(t)$ is the observer gain and $F(t)$ is a feedback gain. State estimation error $e(t)\triangleq x_p(t)-\hat{x}_p(t)$ obeys
\begin{equation}
    \tfrac{d}{dt}e(t) = \left[A_p(t)+L(t)C_p(t)\right]e(t) + B_{pw}(t)w(t).
\end{equation}
In addition, we have that
\begin{equation}
    K(t)x_p(t)-u(t) = \left[K(t)-F(t)C_p(t)\right]e(t)
\end{equation}
and the estimation error covariance $\Sigma_e(t)\triangleq\Ex\{e(t)e^T(t)\}$ is the solution to the Lyapunov differential equation
\begin{multline}
    \tfrac{d}{dt}\Sigma_e(t) = \left[A_p(t)+L(t)C_p(t)\right]\Sigma_e(t) \\ +\Sigma_e(t)\left[A_p(t)+L(t)C_p(t)\right]^T + B_{pw}(t)B_{pw}^T(t) . \label{res_cov_lyap}
\end{multline}
We let $\hat{x}_p(t_i)=0$, and consequently initial condition $\Sigma_e(t_i) = \Theta$. The second term on the right-hand side of \eqref{ltv_perf_obj} is then
\begin{equation}
    \int_{t_i}^{t_f} \textrm{tr}\bigg\{N\left[K-FC_p\right]\Sigma_e \left[K-FC_p\right]^T \bigg\} dt . \label{ltv_perf_obj_2}
\end{equation}

We want to design $L(t)$ and $F(t)$ such that \eqref{ltv_perf_obj_2} (and hence $J$) is minimized, subject to the constraint that the mapping $v\mapsto -u$ is passive. 
Applying Theorem III.3 in \cite{forbes2010passive}, this is guaranteed if there exists a matrix $W(t)=W^T(t) \geq 0,~\forall t\in[t_i,t_f]$ such that
\begin{equation}
    \begin{bmatrix}
    \Pi_1(t) & \Pi_2(t)  \\ \Pi_2^T(t) & \Pi_3(t)
    \end{bmatrix} \leqslant 0 \label{ltv_matrix_ineq}
\end{equation}
where, suppressing time-dependency on all terms, we have
\begin{align}
    \Pi_1 &\triangleq [A_p+LC_p+B_{pu}(K-FC_p)]^T W \nonumber \\ 
    & \quad + W[A_p+LC_p+B_{pu}(K-FC_p)] 
    + \tfrac{d}{dt}W(t), \label{pi1} \\
\Pi_2 & \triangleq W[B_{pu}F-L]+[K-FC_p]^T, \label{pi2} \\
\Pi_3 &\triangleq F+F^T. \label{pi3}
\end{align}

Next, we outline a procedure that can be used to generate feasible gains $L(t)$ and $F(t)$, and assess their performance. 
Suppose that at time $t$, $W(t)$ is given, and $\Pi_1(t),~\Pi_2(t),$ and $\Pi_3(t)$ are specified which satisfy \eqref{ltv_matrix_ineq}. 
Then, \eqref{pi3} is solved for $F(t)$. 
Note that this solution is unique if $F(t)$ is constrained to be symmetric. 
Using $F(t)$, \eqref{pi2} is then solved for $M(t) \triangleq W(t)L(t)$.
Using this solution, \eqref{pi1} is solved for $\tfrac{d}{dt}W(t)$.
In this manner, $W(t)$ may be integrated backward in time, from a final value $W(t_f)=W_f$. 
With trajectories obtained for $W(t),~M(t),$ and $F(t)$,  $L(t)=W^{-1}(t)M(t)$ is computed. 
Using $L(t)$ and $F(t)$, $\Sigma_e(t)$ is found by integrating \eqref{res_cov_lyap} forward in time. 
Finally, performance \eqref{ltv_perf_obj_2} is evaluated. 

The above procedure could be used to optimize over matrix trajectories $\Pi_1(t),~\Pi_2(t)$, $\Pi_3(t)$, and final value $W_f > 0$, for minimal $J$. However, this problem is nonconvex, and intractable due to the very large optimization domain. 
Nonetheless, although still nonconvex, the optimization becomes tractable if $\Pi_1(t),~\Pi_2(t),$ and $\Pi_3(t)$ are constrained to be time-invariant, because the optimization domain is small.
The resulting admittance $\map{Y}:v\mapsto -u$ has state-space realization 
\begin{align}
    A(t) =& A_p(t)+L(t)C_p(t)-B_{pi}(t) C(t) \\
    B(t) =& B_{pi}(t)F(t)-L(t) \\
    C(t) =& F(t)C_p(t)-K(t) \\
    D(t) =& -F(t).
\end{align}

\begin{figure}
  \centering
	\includegraphics[scale=0.55]{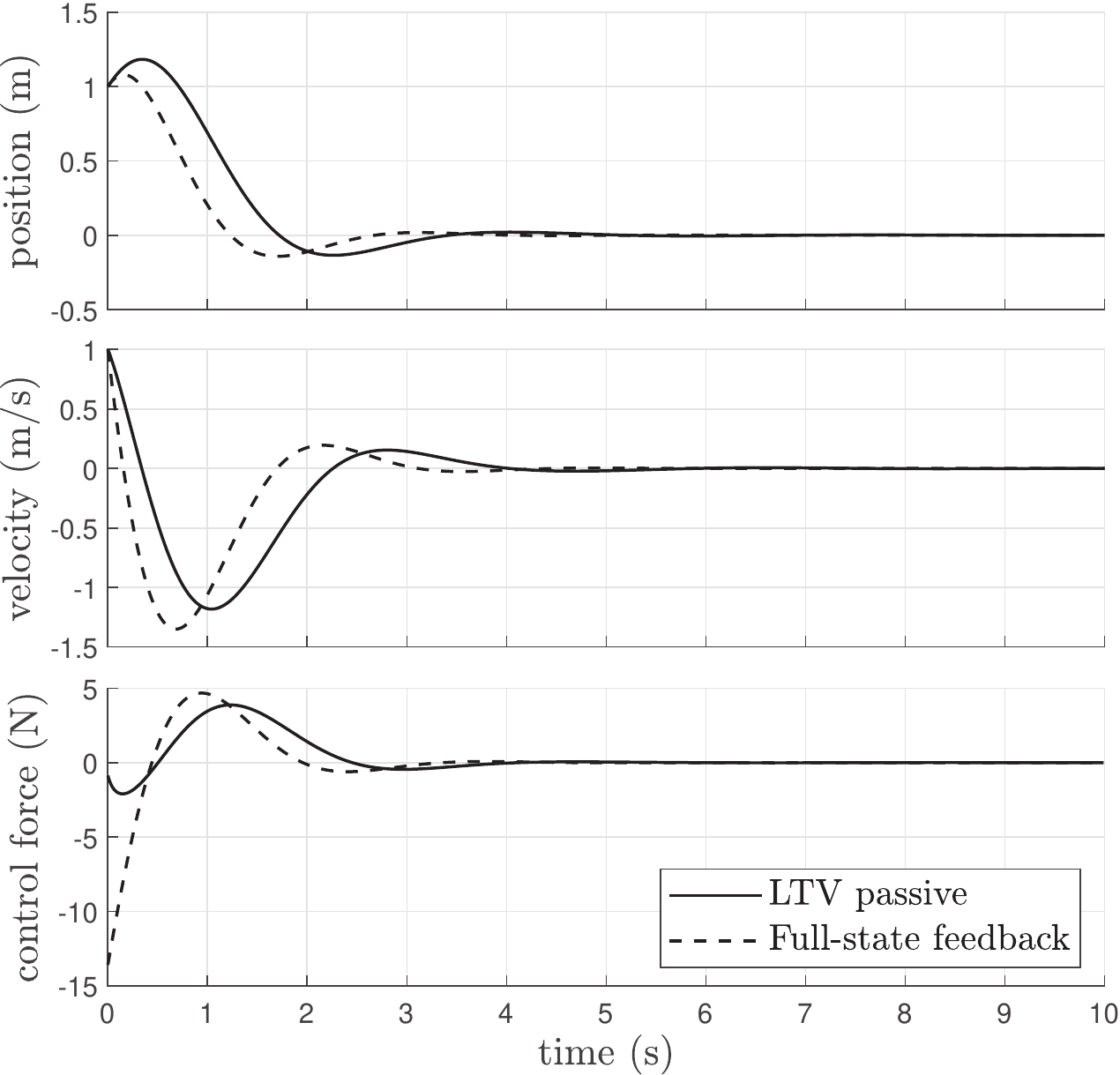}
	\caption{Comparison of time histories for LTV passive and full-state feedback controllers}
	\label{LTV_time_histories}
\end{figure}

Applying this methodology, we assume $m_i=1.5$kg, $m_f=1$kg, $\alpha=0.5$, $c=10^{-5}$Ns/m, $k=5$N/m, $K_e=1$V-s/m, $t_i=0$s, $t_f=10$s, and time step $\Delta t=0.01$s. We fix $S_f = 0$, $W_f = \text{diag}\{1.49\times 10^{-8},1.49\times 10^{-8}\}$, $N(t)=0.75$, $Q(t) = \text{diag}\{100,1\}$, $\Theta = I$, and $B_{pw} = \begin{bmatrix} 0 & 0.1/m(t) \end{bmatrix}^T$.
The nonconvex optimization is performed using MATLAB's {\tt fmincon} function. 
Fig. \ref{LTV_time_histories} shows simulation results, both for the optimized LTV passive controller, as well as for the unconstrained optimal full-state feedback law (i.e., $u(t) = K(t) x_p(t)$). 
Initial conditions $x_{p,1}(0)=1$m and $x_{p,2}(0)=1$m/s were assumed, and zero process noise (i.e., $w(t)=0$). 

\begin{figure}
  \centering
	\includegraphics[scale=0.55]{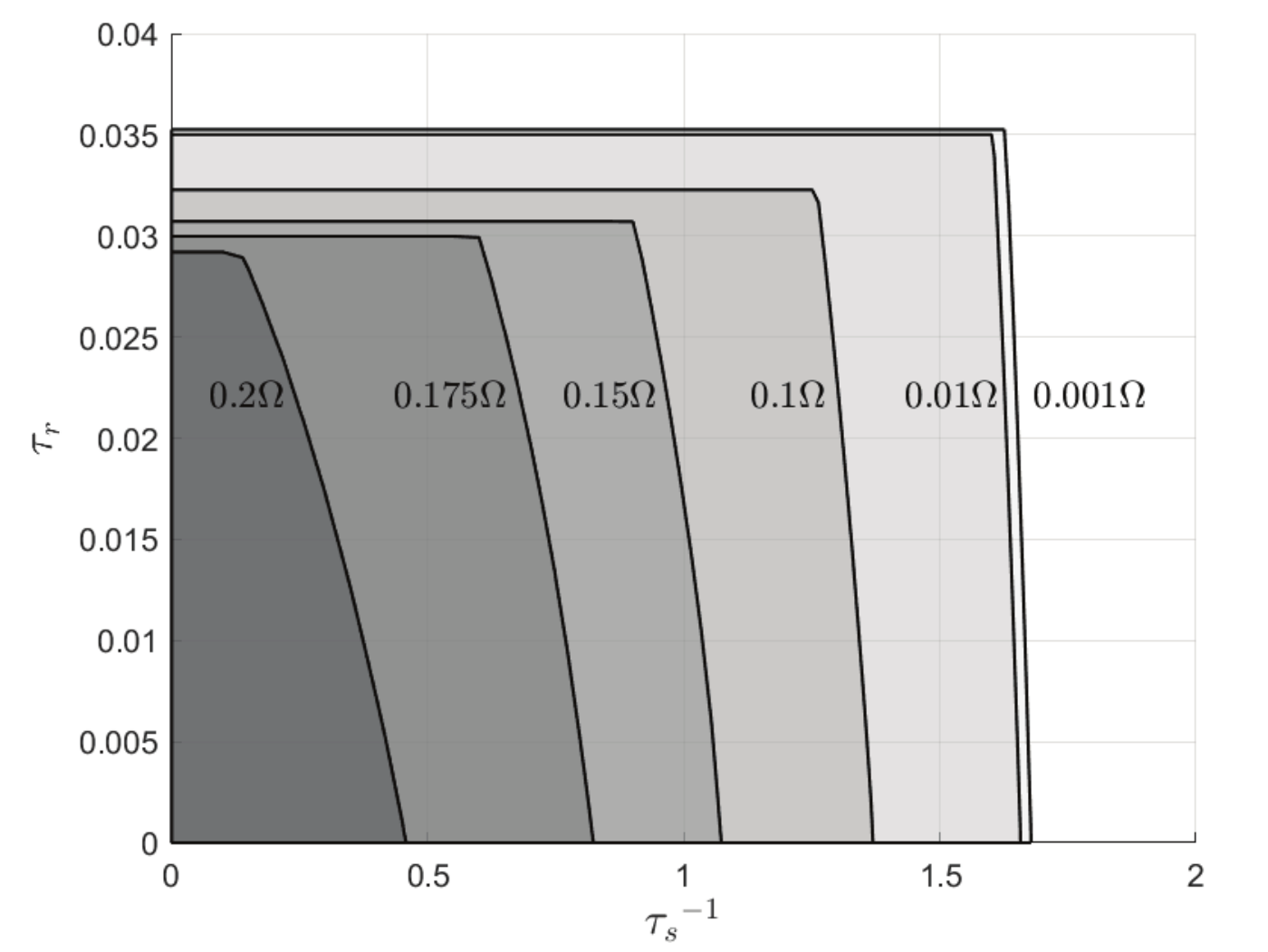}
	\caption{$\{\tau_s^{-1},\tau_r\}$ Pareto fronts (i.e., feasible regions) associated with passive LTV controller $\map{C}$ for various $R$ values}
	\label{LTV_pareto_fronts}
\end{figure}

Now suppose controller $\map{Y}$ is implemented as an SPSA. 
In this case, we are interested in characterizing the $\{R,\tau_r,\tau_s\}$ parameters for which $\map{Y}\in\mathbb{Y}_1(R,\tau_r,\tau_s)$. 
Technically, feasibility would require the existence of continuous matrix functions $P(t) = P^T(t) > 0$ and $X(t)$ such that \eqref{ltv-finite-dim-theorem-lmi1} and \eqref{ltv-finite-dim-theorem-lmi2} hold for all $t \in [t_i,t_f]$. 
Here, we approximate this feasibility criterion by enforcing it at a large-but-finite set of evenly-spaced discrete times on this interval. 
The optimization domain consists of the values of $P(t)$ and $X(t)$, evaluated at these discrete times. 
Finite-difference approximation is then used to evaluate $\tfrac{d}{dt}P(t)$ at the discrete times. 
As such, we have the following optimization problem:
\begin{equation*} 
\begin{array}{lll}
\text{Given} &:& R,\tau_s, A_k,B_k,C_k,D_k \text{ for } k \in \{0,...,N_s\} \\
\text{Maximize} &:& \tau_r \\
\text{Over} &:& \tau_r,P_k,X_k, \text{ for } k \in \{0,...,N_s\} \\
\text{Subject to} &:& \eqref{ltv-finite-dim-theorem-lmi1},\eqref{ltv-finite-dim-theorem-lmi2}, P_k=P_k^T>0
\end{array}
\end{equation*}
Cross sections of the resultant Pareto surface for $(R,\tau_s,\tau_r)$ are shown in Fig. \ref{LTV_pareto_fronts}. Note the general trends are very similar to Fig. \ref{LTI_pareto_fronts}.
As $\R$ is increased, the feasible region of $\{\tauR,\tauS^{-1}\}$ values shrinks, implying that as the actuator is made less efficient, the storage system must be made more efficient.

\subsection{Fractional-order SPSAs}

In this final example, we explore the self-powered feasibility of infinite-dimensional admittances. Specifically, we consider the SISO fractional-order (FO) lead-lag filter \cite{dastjerdi2019linear} characterized by transfer function
\begin{equation}
    \hat{Y}(s) = k_p \left( \frac{1+s/\omega_L}{1+s/\omega_h}\right)^\mu \label{FO_eqn}
\end{equation}
where $\omega_L < \omega_h$ and $\mu \in (-\infty,\infty)$. Note $\hat{Y}(s)$ is infinite dimensional for non-integer $\mu$.

We first reformulate $\hat{Y}(s)$ as the LFT in Fig.~\ref{lft} where 
\begin{equation} \label{inf_dim_Z}
    Z(t) = \begin{bmatrix}
    z_{11} & z_{12} \\
    z_{21} & z_{22}
    \end{bmatrix}
\end{equation}
and
\begin{equation}
\hat{G}(s) = -\frac{\hat{Y}(s)-z_{11}}{z_{12}z_{21}+z_{22} (\hat{Y}(s)-z_{11})}
\end{equation}
where $z_{12},z_{21},z_{22}$ are time-invariant scalar parameters, subject to constraint \eqref{Z(t)-norm-condn}, and 
\begin{equation} \label{z11_eqn}
    z_{11} = k_p \left(\omega_h/\omega_L\right)^\mu.
\end{equation}

\begin{figure}
  \centering
	\includegraphics[scale=0.55]{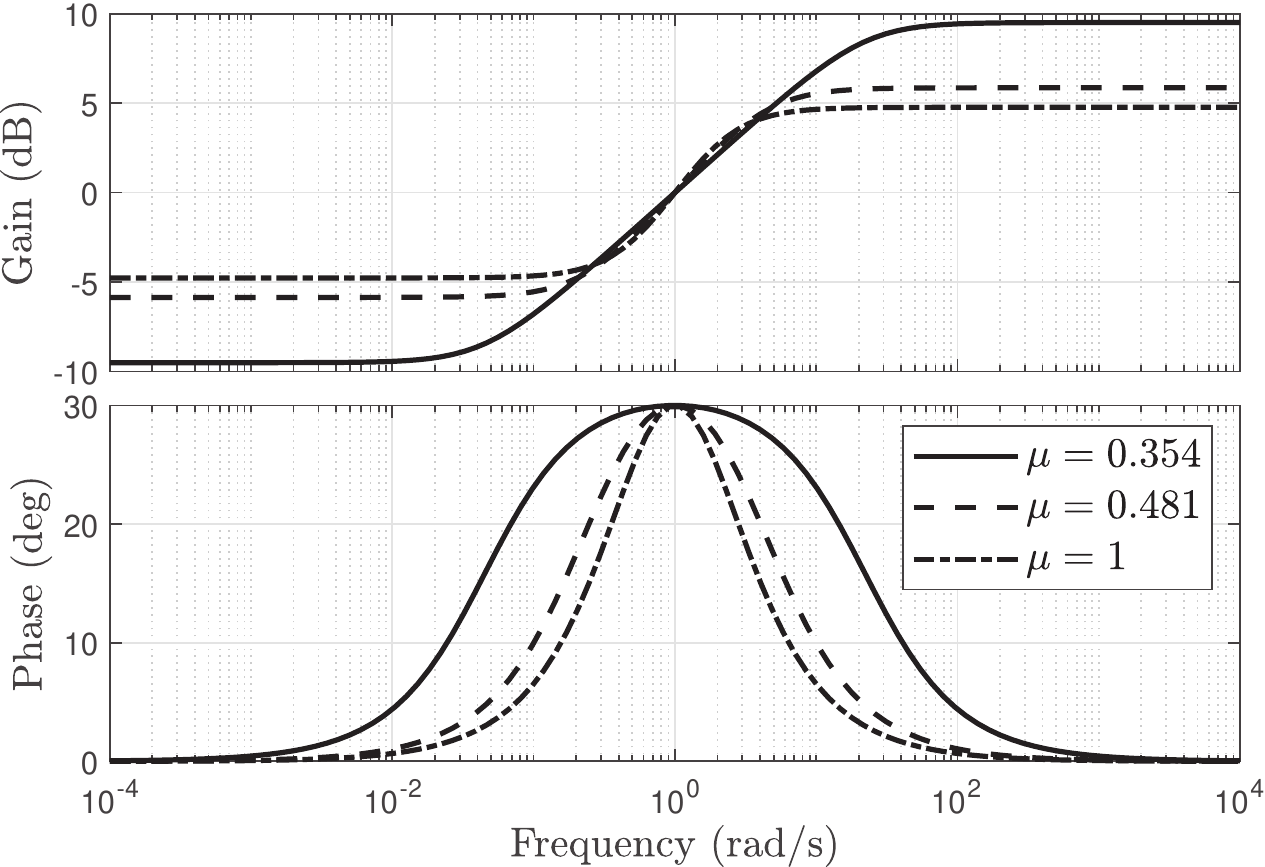}
	\caption{Achieving peak phase and gain specifications with varying $\mu$}
	\label{mu_effect}
\end{figure}

Suppose we would like to realize FO lead-lag filters that attain a certain peak phase lead (or lag) $\phi_m$ and gain $\gamma_m$ at a particular frequency $\omega_m$. It is possible to achieve such specifications using different combinations of the filter parameters $\{\mu,k_p,\omega_L,\omega_h\}$. In fact, for given $\{\phi_m,\gamma_m,\omega_m,\mu\}$ we can uniquely compute the necessary values of $\{k_p,\omega_L,\omega_h\}$.  For example, suppose we specify that $\phi_m = +30^\circ$ , $\gamma_m = 1$, and $\omega_m = 1$ rad/s. Fig. \ref{mu_effect} shows the resulting bode plots of FO lead filters corresponding to various $\mu$ values. As $\mu$ decreases, there is a reduction in the DC gain, an increase in high-frequency gain, and a widening of the phase lead.

In the context of SPSAs, it is natural to ask how varying $\mu$ in this way affects the least efficient energy parameters (i.e., the largest combinations of $R,\tau_r,\tau_s^{-1}$) required to realize the corresponding filter. To answer this question, we begin by observing that \eqref{inf-dim-lti-G-condition} is equivalent to 
\begin{equation}\label{new_pr_constraint}
    z_{12}z_{21}\mathbf{Re}\{\hat{\tilde{Y}}(j\omega)-z_{11}\} + z_{22}|\hat{\tilde{Y}}(j\omega)-z_{11}|^2 \leqslant 0, ~\forall \omega\in\mathbb{R}
\end{equation}
where $\hat{\tilde{Y}}(j\omega)\triangleq\hat{Y}(j\omega-\tau_s^{-1})$ and $\mathbf{Re}\{\cdot\}$ denotes the real part of a complex number. Furthermore, via application of the initial value theorem, \eqref{inf-dim-lti-G0-condition} is equivalent to
\begin{equation}\label{new_g0_constraint}
    \frac{\tau_r z_{11}\mu\left(\omega_h-\omega_L\right)}{z_{12}z_{21}} \leq 1.
\end{equation}

In a manner similar to the two previous examples, we now formulate an optimization problem that maximizes $\tau_r$ given energy storage subsystem parameters $\{R,\tau_s\}$. For this example, the corresponding optimization is:
\begin{equation} 
\begin{array}{lll}
\text{Given} &:& k_p,\omega_L,\omega_h,\mu,R,\tau_s \\
\text{Maximize} &:& \tau_r \\
\text{Over} &:& \tau_r,z_{12},z_{21},z_{22} \\
\text{Subject to} &:& \eqref{Z(t)-norm-condn}, \eqref{FO_eqn},\eqref{inf_dim_Z},\eqref{z11_eqn},\eqref{new_pr_constraint} ,\eqref{new_g0_constraint}
\end{array}
\end{equation}
Although not immediately apparent, this optimization is convex. We note that constraint \eqref{new_pr_constraint} is enforced at a finite number of evenly-spaced frequencies. This frequency vector should obviously be of sufficient length and resolution to accurately capture the dynamics of $\hat{Y}(s)$.

\begin{figure}
  \centering
	\includegraphics[scale=0.55]{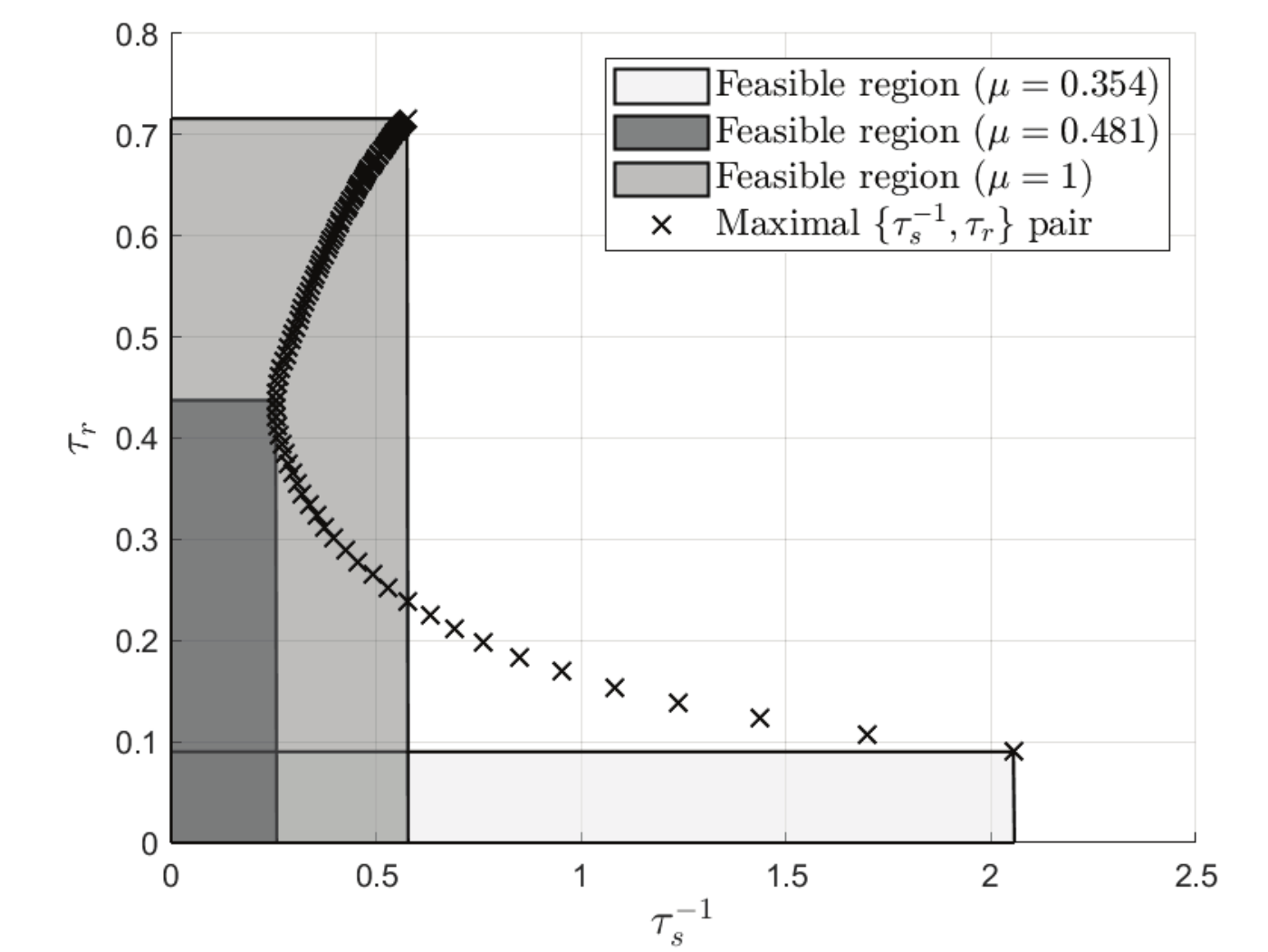}
	\caption{Effect of $\mu$ on maximal $\{\tau_s^{-1},\tau_r\}$ parameters for $R=0.1\Omega$}
	\label{backbones}
\end{figure}

\begin{figure}
  \centering
	\includegraphics[scale=0.55]{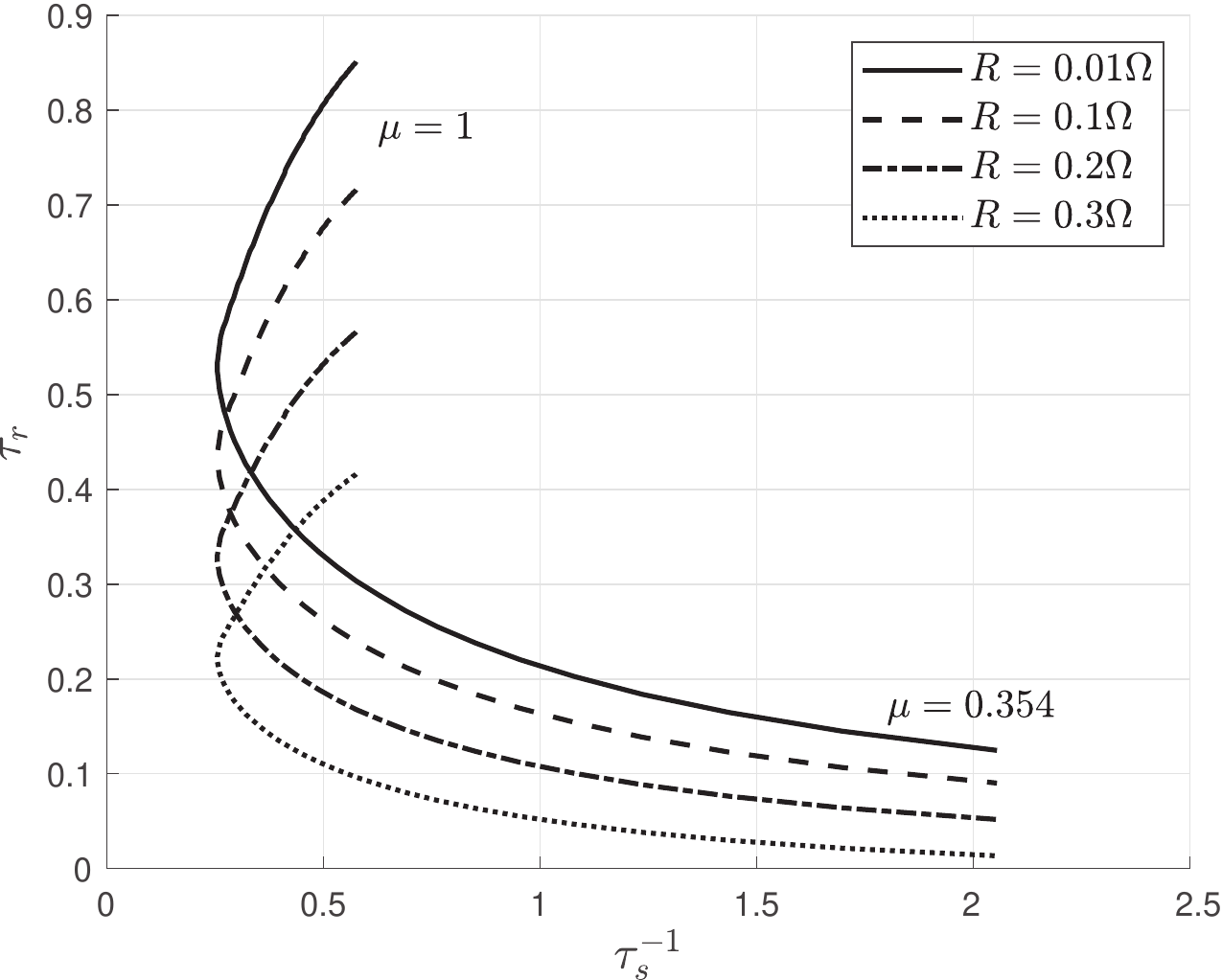}
	\caption{Effect of $\mu$ and $R$ on maximal $\{\tau_s^{-1},\tau_r\}$ parameters}
	\label{backbones2}
\end{figure}

Returning to our initial question, we let $R=0.1\Omega$ and construct the $\{\tau_s^{-1},\tau_r\}$ Pareto fronts associated with values of $\mu$ ranging from $0.354$ to $1$ and the corresponding $\{k_p,\omega_L,\omega_h\}$ parameters necessary to achieve $30^\circ$ phase lead and unity gain at $1$ rad/s. Three typical feasible regions (corresponding to the $\mu$ values used in Fig. \ref{mu_effect}) are shown in Fig. \ref{backbones}. Note that these regions are rectangular, implying that the maximal $\tau_r$ and $\tau_s^{-1}$ parameters are independent for a given $\mu$. Plotting just the corner points of these regions yields a smooth ``backbone" curve, also shown in Fig. \ref{backbones}. We repeat this process for increasing values of $R$ and compare the resulting backbones in Fig. \ref{backbones2}. The downward shift of the  backbones suggests that $\tau_r$ is coupled with $R$ but $\tau_s$ is not.

\section{Conclusions} \label{sec:sec5}

Self-powered feedback control is a unique and attractive approach for controlling a variety of physical systems. In this paper, we derived sufficient conditions for a synthetic admittance to be self-powered feasible given certain parasitic loss parameters. These conditions are broadly applicable to infinite-dimensional and time-varying linear admittances. Furthermore, three numerical examples were presented, which demonstrated how the feasibility conditions may be used to identify a Pareto front of the least-efficient loss parameters necessary to realize a given admittance as an SPSA. 

There is much future work to be done in the area of self-powered systems. One obvious next step is the identification of feasibility conditions which are both necessary and sufficient. In addition, the theory developed herein should be extended to accommodate more realistic parasitic loss models. Lastly, as our focus in this paper was restricted to linear feedback control laws, research on the feasibility and synthesis of nonlinear self-powered controllers is also warranted.

\bibliographystyle{IEEEtran}
\bibliography{regen}

\end{document}